\documentclass[11pt]{amsart} \usepackage{amsfonts} \usepackage{amscd}
\usepackage{amssymb}
\newcommand{\ring}[1]{\mathbb{#1}}
\newcommand{\C}{\ring{C}} \newcommand{\Q}{\ring{Q}}
\newcommand{\Z}{\ring{Z}} 
 \newcommand{\G}{\ring{G}}
\newcommand{\bu}{\bullet}
\newcommand{\ms}{\medskip}

\newcommand{\ra}{\rightarrow}

\newcommand{\iso}{\xrightarrow{\sim}}
\newcommand{\xra}{\xrightarrow}
\newcommand{\gm}{\G_m}
\newcommand{\nd}{\noindent}

\newcommand{\ub}{\underline}
\newcommand{\mh}{\mathcal M\mathcal H_{pol}^1}
\newcommand{\hra}{\hookrightarrow}

\def\1{{\mu\mkern-6mu\mu}}

\title[one-motives] 
{Duality of Albanese and Picard 1-motives} 
\author{Niranjan Ramachandran}
\address{Department of Mathematics, University of Maryland, College
Park, MD 20742-4015, USA} 
\email{atma@math.umd.edu}
\date{30 September 2000}

\begin{document}
\newtheorem{theorem}{Theorem}[section]
\newtheorem{lemma}[theorem]{Lemma}
\newtheorem{corollary}[theorem]{Corollary}
\newtheorem{proposition}[theorem]{Proposition}
\theoremstyle{definition}
\newtheorem{remark}[theorem]{Remark}
\newtheorem{definition}[theorem]{Definition}
\newtheorem{example}[theorem]{Example}

\begin{abstract}
 We define Albanese and Picard 
1-motives of smooth (simplicial) schemes over a perfect field. For smooth
proper schemes, 
these are the classical  Albanese and Picard varieties. For a curve, these are the homological 1-motive
 of Lichtenbaum and 
the motivic $H^1$ of Deligne. This paper proves a
conjecture of 
Deligne about providing an algebraic description, via 1-motives, 
of  the first homology and 
cohomology groups of a complex algebraic variety. (L. Barbieri-Viale
and V. Srinivas have also proved this independently.) It also contains a purely
algebraic proof of Lichtenbaum's conjecture that the Albanese and  
the Picard 1-motives of a (simplicial) scheme are dual. This gives a
new proof of an unpublished theorem of Lichtenbaum that Deligne's
1-motive of a curve is dual to Lichtenbaum's 1-motive.
\end{abstract}
\thanks{Supported by NSF and the Alfred P. Sloan Doctoral Dissertation
Fellowship, 1995-96}
\maketitle

\setcounter{section}{-1}

{\bf Revised version}

\section{Introduction}

Let $k$ be a perfect field. Two classical invariants of any smooth
projective variety $X$ over $k$ are the abelian varieties $Alb(X)$
(the Albanese variety), and $Pic(X)$ (the Picard variety) \cite{co}. It is a
classical observation of W. L. Chow, J. Igusa, and A. Weil that
$Alb(X)$ and $Pic(X)$ are naturally dual. 

In this paper, we extend this classical result, using 1-motives,  to arbitrary schemes
over $k$.  P. Deligne's technique  \cite{h} \S 6.3 of
defining a cohomological invariant for arbitrary varieties consists of
defining it for simplicial schemes $X_{\bu}$ (this is the complement
of a simplicial divisor with normal crossings $E_{\bu}$ in a smooth
projective simplicial scheme $X'_{\bu}$) and using smooth
hypercoverings to define it for varieties (see \S \ref{positive}). For
any such scheme $X_{\bu}$,  we
will define 
the Albanese 1-motive $M_1(X_{\bu})$ and the Picard 1-motive
$M^1(X_{\bu})$. Our construction of the covariant functor $M_1$ is
based on the generalized Albanese variety of
Rosenlicht-Lang-Serre. The construction of the contravariant functor
$M^1$  is based upon the Picard scheme of a simplicial scheme.

Our main result (Theorem \ref{duality}) is that, for any $X_{\bu}$, the 1-motives $M^1$ and
$M_1$ are naturally dual. This was conjectured by S. Lichtenbaum. When $X_{\bu}= X'_{\bu}$ is a constant scheme corresponding to a smooth projective
variety $V$ (i.e., each $X_i = V$), Theorem \ref{duality} reduces to the
classical duality of $Alb(V)$ and $Pic(V)$. 

For
any curve $C$, Deligne has defined a
1-motive $H^1_m(C)(1)$  and obtained, via the realization
functors on one-motives, the $\ell$-adic cohomology, De
Rham cohomology etc., of $C$ from $H^1_m(C)(1)$. Lichtenbaum \cite{L1}
showed that a 1-motive appears,  rather
naturally, in the Suslin complex of the curve $C$. This led him to
define $h_1(C)$, the homological 1-motive of $C$. The homological
1-motive of a smooth curve is readily seen to be a generalized Jacobian 
 of M. Rosenlicht \cite{Ro1, Ro2}. 

Associated with any curve $C$ is a canonical smooth simplicial scheme
$Z_{\bu}$ (\S \ref{kurv}; \cite{h} \S 10.3). We show that $M_1(Z_{\bu}) = h_1(C)$ and
$M^1(Z_{\bu}) = H^1_m(C)(1)$. As a consequence of Theorem \ref{duality}, we
obtain a new proof of an unpublished  theorem of Lichtenbaum (Theorem
\ref{slick}): \emph{the 1-motive $h_1(C)$
is naturally dual to the 1-motive $H^1_m(C)(1)$.}  

Let $X_{\bu}$ be a scheme (as above) over $\C$. We also  show that the Hodge realization of $M^1(X_{\bu})$  is the  
mixed Hodge structure $H^1(X_{\bu};\Z(1))$. 
It follows from the duality theorem (Theorem \ref{duality}) that the
Hodge realization of $M_1(X_{\bu})$  is  the mixed Hodge structure on $H_1(X_{\bu})/{torsion}$.
This settles the conjecture of Deligne for $I(H^1(X_{\bu})(1))$ (\S
10.4.1  of
\cite{h}). Results about the other realizations of these
1-motives will appear elsewhere.

Weak Lefschetz-type theorems are available for varieties which are
either non-singular quasi-projective or affine (\cite{gm} pp. 23-24);
it is also remarked in (loc. cit) that such theorems do not hold in
general for arbitrary
varieties. Therefore, the study of the motivic $H^1$ and
$H_1$ for general varieties does not seem to reduce to the case of
curves. This means that the results obtained  here about the motivic $H^1$ and
$H_1$ are perhaps of interest. Related work in this direction includes
 \cite{on}, \cite{FW}, \cite{bs}, \cite{bs2}, \cite{esv}, \cite{li3}, and \cite{me2}. I
 should mention that \cite{me2} generalizes to $H^n$ (for $n >1$)
\footnote{In \cite{me2}, the conjecture of Deligne is proved up to isogeny for
$H^n(V)$; here, our theorem proves the integral version of the
conjecture for $H^1$.} the results 
 obtained here for $H^1$. But the Albanese
 1-motive and the duality theorem are not touched there. Recently
L. Barbieri-Viale and 
 V. Srinivas \cite{bs2} (announced in \cite{bs}) have independently
proved Theorems \ref{repic} and 
\ref{picho}. One (the essential) difference
between their work and ours is that their Albanese 1-motive is
\emph{defined} to be the dual of the Picard 1-motive. Another is that they define
their invariants only in characteristic zero, whereas we do so for
an arbitrary perfect field (The restriction to perfect fields is
forced by Theorem \ref{reddy}).

\subsection{Notations}
 We use the following conventions.\ms

{\small
$S = $ Spec~$k$, the spectrum of  a perfect field $k$.

The characteristic of $k$ is denoted by $p$.

$\bar{S}$ = Spec~$\bar{k}$, the spectrum of an algebraic closure
$\bar{k}$ of $k$.

{\bf g} = $Gal(\bar{k}/k)$, the Galois group of $\bar{k}$ over $k$.

 We consider only separated schemes $X$ over $S$.  Such schemes are supposed to
be of finite type over $S$ unless otherwise mentioned.  We use
$\bar{X}$ to denote the $\bar{S}$-scheme $X \times_S \bar{S}$ defined
 by the
$S$-scheme $X$. 

$\tilde{X}:=$ the normalization of $X$.

$X^{red}:=$ the maximal open and closed (reduced) subscheme of $X$.

$w(X)$:= the set of connected components of $\tilde{X}\times_S \bar{S}$ 

 = the set of  irreducible components of $\bar{X}$.

 A variety is a scheme $V$ over $S$ such that $\bar{V}$
is a integral scheme.

We use $S$ and $k$ interchangeably: we say $k$-points of a scheme to
mean its
points with values in $S$.

The group schemes we consider are commutative. For such a group scheme
$G$, let $G^0$ denote the connected component of the identity and
$G^{0, red}$ denote the connected component of the identity of the
reduced scheme 
$G^{red}$.  
For a positive integer $r$ prime to $p$, we denote by $_r G$ the
 group scheme corresponding to
$Ker( G \xra{r} G)$,
 i.e., the $r$-torsion subgroup scheme of $G$. 

A semiabelian scheme is a commutative group scheme
which is an extension of an abelian scheme by a torus.

$\ub{\Z}$ is the locally algebraic group scheme defined
by $\ub{\Z} (\bar{S}) = \Z$ with trivial {\bf g} action.

$\mathbb G'$~~ := the category of group schemes over $S$.

$\mathbb G$~~ := the category of locally algebraic
semiabelian schemes over $S$.

$\mathbb D$ := the abelian category of sheaves (of abelian
groups) on $S_{FL}$.

$S_{FL}$ := the (big) flat site of $S$ \cite{Mi}.

Properties of a simplicial scheme $X_{\bu}$ are properties that hold for each 
$X_i$: $X_{\bu}$ is proper if each $X_i$ is proper over $S$, same for smooth, 
normal, reduced, etc.

Let $W$ be a set, and denote by $\Z[W]$ the free abelian group
generated by $W$.  We define the abelian group $\Z^W$ to be the group of integer valued 
functions on $W$ i.e. $\Z^W:= Maps(W,\Z)$. For a finite set $W$, the groups 
$\Z(W)$ and $\Z^W$ are naturally dual as abelian groups.

 For any finitely generated abelian group $G$, we put $G_{tor}$ to denote the
 torsion subgroup of $G$. We denote by $G/{tor}$ the maximal torsion-free
 quotient of $G$.
  
  $\mathcal M \mathcal H \mathcal S$:= the abelian category of mixed
  Hodge structures (MHS). 

$\mh$ := the category of torsion-free MHS $N$ of type $$\{ (0,0), (-1,0),
(0,-1), (-1, -1)\}$$ such that $Gr^W_{-1}N$ is polarizable \cite{h} \S
10.1. 

$\mathfrak T(M) \in \mh$ denotes the Hodge realization (\cite{h} 10.1)
of a 1-motive
$M$.  
  
We use $H_{\Z}$ to denote the underlying
$\Z$-lattice of the mixed Hodge structure $H$. We use $W_iM$ and
$W_iH$ to indicate the weight filtration on 1-motives and mixed
Hodge structures. The Tate Hodge
structure $\Z(1)$ corresponds to the lattice $2\pi i \Z \subset
\C$; the Hodge structure $\Z(n)$ (for $n>0$) is defined to be
$\Z(1) ^{\otimes n}$.  For $n<0$, we define $\Z(n)$ to be
the Hodge structure dual to $\Z(-n)$ i.e we have $\Z(n) =
Hom_{MHS}(\Z(-n), \Z)$. The Tate Hodge structure
$\Z(n)$ is of weight $-2n$ of type $(-n, -n)$. For any mixed
Hodge structure $H$, the Tate twist $H(n)$ denotes   $H \otimes\mathbb Z(n)$.
For any scheme $V$ over $\C$, we use $H^*(V)$ to indicate the
mixed Hodge structure $H^*(V(\C);\Z)$. Similar convention for
homology applies as well.  Given  a sheaf $\mathcal F$ on $V$, we
write $H^i(V;\mathcal F)$ for $H^i(V(\C);\mathcal F)$.}\ms

\noindent {\bf Acknowledgements.} {\small It is a pleasure for me to
  thank S. Lichtenbaum for much 
encouragement. The results 
presented here would not 
have been possible without his guidance; that they appear in published form is
due to his insistence. This paper is a greatly
revised version of \cite{R1}.  I would like to thank the mathematics
department at Brown, the Alfred P. Sloan foundation, 
and the National Science Foundation (NSF) for their  
support. Lastly, I heartily thank the referee for the careful reading,
enthusiasm, and numerous suggestions thereby improving the exposition.}\medskip

\section{Albanese Schemes}

This section is based on a succinct suggestion of J.-P. Serre
\cite{se4}. The goal is the definition of the Albanese scheme of a
scheme $X$ over $S$. This is different from the (very interesting)
universal regular 
quotient of \cite{esv}.

\subsection{Locally algebraic group schemes.} 
\begin{definition} A locally algebraic group scheme is a
group scheme locally of finite type over $S$. \end{definition}  
 
A (group) scheme of finite type will sometimes be called an algebraic
(group) scheme. 
   A locally algebraic group scheme $G$ such that $G^0 = S$ is said to be
discrete or of {\it dimension zero}; it is an \'etale $S$-scheme. Such a
group $G$ is determined 
(up to isomorphism) by the discrete {\bf g}-module $G(\bar{S})$.  

Recall from \cite{DG} pp. 284-6 that there is a discrete locally algebraic group $\pi_0(G)$ associated with
a locally algebraic group scheme $G$, a scheme-theoretic analog of the
group of connected components; it is the universal object for
homomorphisms from $G$ to discrete locally algebraic group schemes.
One has an exact sequence of locally algebraic group schemes:
\begin{equation}\label{qg}
0 \ra
G^0 \ra G \xra{q_G} \pi_0(G) \ra 0 \end{equation} where $q_G$ is a
homomorphism of group schemes.

\subsection{Construction}\label{bps}
 Let $B$ be a commutative group scheme (algebraic) over
 $S$. Given two $B$-torsors $P_1$ 
and $P_2$ (torsor = principal homogeneous space),
  the $B$-torsor obtained by Baer summation is
denoted $P_1 \vee_B P_2$. For any $B$-torsor $P$, one can define a new
 torsor $P^{-1}$ by twisting the action of $B$ by the automorphism
 corresponding to the inverse $i:B \ra B$.
 With each $B$-torsor $P$, 
 one can associate a locally
 algebraic group $B_P$ as follows. For $n$
positive, let $P^{\otimes n}$ denote the $B$-torsor $P \vee_B
\cdots\vee_B P$ ($n$ terms).  For $n$ negative,  let
$P^{\otimes n}$ denote the $B$-torsor $P^{-1} \vee_B \cdots \vee_B P^{-1}$
($n$ terms). Let $P^{\otimes 0}$ be the trivial torsor $B$. We define   
\begin{equation}
 B_P:= \coprod_{n \in \Z} P^{\otimes n}; \end{equation}
note the canonical morphism $a_P:B_P \twoheadrightarrow \ub{\Z}$. The natural morphisms
\begin{equation}
P^{\otimes n}\times_S P^{\otimes m} \ra P^{\otimes n}\vee_B P^{\otimes
 m}:=  P^{\otimes (m+n)} 
\end{equation}
 endow $B_P$ with the structure of a commutative 
group scheme; $B_P$ is an extension of $\ub{\Z}$ by the group 
scheme $B$. Any scheme $T$ defines a sheaf $\mathcal T$   of sets on
$S_{FL}$:  $X \mapsto {\rm Hom}_S(X,T)$.  There 
 is a natural inclusion of $\mathcal P \hra \mathcal B_P$ of sheaves
 of sets on 
$S_{FL}$; in fact, $\mathcal P$ is the subsheaf of $\mathcal B_P$ 
corresponding to the inverse image under $a_P$ of the section $1\in
\ub{\Z}$.    

 A semiabelian scheme is an
algebraic group scheme which is an extension of an abelian scheme by a
torus; it is geometrically connected. We call a locally algebraic
group scheme $G$ {\it semiabelian} if 
$G^0$ is a semiabelian scheme. If $B$ is a semiabelian scheme
and $P$ is a torsor under $B$, then the group scheme $B_P$ is a
semiabelian locally algebraic  group scheme.\ms

\begin{definition} A large group is a semiabelian locally algebraic
 group  scheme. \end{definition} 
     
We use this terminology because (1) the name is short and (2) a large group
is not a scheme of finite type. We denote by $\mathbb D$ the abelian category of sheaves (of abelian
groups) on the (big)
flat site $S_{FL}$ of $S$. The category $\mathbb G$ of large groups
embeds as  a
full  subcategory (denoted $\mathbb G$) of $\mathbb D$. 

\subsection{Morphisms to large group schemes} 
 For any set
$W$, let $\Z[W]$ be the free abelian group generated by the
set $W$. 
Let $X$ be a  scheme.  Define $P_X$ to be the presheaf of abelian groups on
$S_{FL}$ by $T \mapsto \Z[X(T)]$.  

\begin{definition}  $\mathcal Z_X$ is the sheaf on $S_{FL}$ associated
with the presheaf $P_X$. \end{definition}

The association $X \mapsto \mathcal Z_X$ determines a covariant functor
from the category of schemes over $S$ to $\mathbb D$.
If we consider $X$ as a sheaf $\mathcal X$ of sets, then we have a canonical inclusion $j: \mathcal X \hra
\mathcal Z_X$.  The sheaf $\mathcal Z_X$ has the following
universal property:  
Given any abelian sheaf $\mathcal F$ on $S_{FL}$ and a map $g: \mathcal X \ra
\mathcal F$ of sheaves of sets, there is a unique homomorphism of abelian
sheaves $\tilde{g}: \mathcal Z_X \ra \mathcal F$ which factors $g$,
i.e., $\tilde{g}j~=~g$.

It is known that the sheaf $\mathcal Z_X$ is not 
representable in any reasonable topology \cite{voe}.
  One is naturally led to consider the universal representable
quotient of $\mathcal Z_X$; we restrict ourselves to the situation where the
representing scheme is a large group. 

It is well known that there is an \'etale algebraic $S$-scheme
$\pi_0(X)$ and a $S$-morphism $q_X: X \ra \pi_0(X)$ which is universal for
morphisms from $X$ into \'etale $S$-schemes (\cite{DG}, Prop 6.5,
p.154). 

  Consider the sheaf $\mathcal Z_{\pi_0(X)}$ on $S_{FL}$ defined by the
scheme $\pi_0(X)$; as above, $\tilde{q}_X$ denotes the induced
homomorphism from $\mathcal Z_X \ra \mathcal Z_{\pi_0(X)}$.  The following
lemma is well known (loc. cit).

\begin{lemma}\ms

{\rm (a)}  $\mathcal Z_{\pi_0(X)}$ is represented by a
discrete large group $D_X$.

{\rm (b)} The pair $(\mathcal Z_{\pi_0(X)}, \tilde{q}_X)$ is
universal for morphisms from $\mathcal Z_X$ to discrete group schemes.
\end{lemma}
 The group $D_X$ of a variety $X$ is isomorphic to $\ub{\Z}$.


\begin{definition}\label{alba} The \emph{Albanese sheaf} of $X$ is
the universal object for morphisms from $\mathcal Z_X$ to objects in the
subcategory $\mathbb G$ (i.e. sheaves corresponding to large groups).
\end{definition}

We say that $X$ admits a {\it universal morphism} if the {\it Albanese
sheaf} of $X$ exists. 

Grant, for the moment, the existence of the Albanese sheaf of $X$. The
sheaf $\mathcal A_X$ comes equipped with a {\it surjective} homomorphism 
$u_X: \mathcal Z_X \ra \mathcal
A_X$. Let $A_X$ be the large group corresponding to the sheaf $\mathcal
A_X$. The morphism $\mathcal X \hra \mathcal Z_X \ra \mathcal A_X$ of
sheafs of sets on 
$S_{FL}$ corresponds, by Yoneda's lemma, to a morphism $u_X: X \ra
A_X$ of schemes. 

We call $A_X$ the {\it Albanese
scheme} of $X$ and $u_X$ is the universal morphism associated with
$X$. The pair ($A_X, u_X$) is uniquely determined (up 
to unique isomorphism).  An easy consequence of definition \ref{alba} is a
canonical isomorphism \begin{equation} \pi_0(A_X) \iso D_X. {}
\end{equation}  The
assignment $X \mapsto A_X$ determines a covariant functor $A$ from the
category of  algebraic schemes over $S$ to large groups over $S$.

 For a  scheme $W$ finite over $S$ (e.g., a
 finite group scheme), we 
 have canonical isomorphisms
\begin{equation}\label{isos}
\mathcal Z_W \iso A_W \iso D_W. 
 \end{equation}
The next proposition follows  from
 the corresponding properties of the functor $P \mapsto \mathcal Z_{P}$. 
\begin{proposition}\label{enjoy} {\rm (Properties of the
 Albanese scheme)}

{\rm (i) (base change)} {\it Let $S'$ be the spectrum of a perfect
  field $k'$ containing $k$. One has $A_X \times_S S' = A_{X \times_S S'}$.}

{\rm (ii)} {\it Let $X:= Y \coprod V$ be the disjoint union of
schemes $Y$ and  $V$.  One has $ A_X = A_Y \times_S A_V$.}   

{\rm (iii)} {\it Let $X \xra{f} Y$ be a surjective morphism. The induced  map
$f_*: A_X \ra A_Y$ is a surjection.}
\end{proposition} 

\begin{remark}\label{serre} Serre has shown that every scheme $X$
  admits a universal morphism to a 
torsor $P$ under a semiabelian scheme $B$:  in \cite{Se2}, he shows that every
scheme $X$ admits a universal morphism (over $\bar{k}$) to a semiabelian scheme $B$ (the scheme $B$ is defined
over $k$) and the method of descent as in \cite{Se1} (pp.102-107)
assure the existence
of a universal map (defined over $k$) to a torsor $P$ under the
universal semiabelian scheme $B$. This semiabelian scheme $B$ is the
``generalized Albanese'' of $X$.\end{remark}  

\begin{remark}\label{torso}
Let $T$ (resp. $R$) be a torsor under a group scheme $G$ (resp. under
$H$). Consider the group schemes $G_T$ and $H_R$ together with
the canonical 
surjective homomorphisms $a_T: G_T \ra \ub{\Z}$ and $a_R: H_R
\ra \ub{\Z}$.  
Specifying a  morphism $h$ of torsors between $T$ and $R$ is  
equivalent to specifying a homomorphism of group schemes $\tilde{h}: G_T \ra
H_R$ such that $a_R \tilde{h} = a_T$ (i.e. $\tilde{h}$ induces
the identity on the quotient $\ub{\Z}$).\end{remark}

Let $V$ be a connected scheme (so that there is a map $g_V: D_V \to
\ub{\Z}$).     
Let $f: V \ra P$ be the universal morphism from $V$ to torsors under
semiabelian schemes (as in Remark \ref{serre}). Consider the natural inclusion
of $P$ in $B_P$. The map $f$ induces a map $\tilde{f}: \mathcal Z_V \ra
B_P$.  Define $B'':= B_P\times_{\ub{\Z}} D_V$, where the
fiber product is taken via the maps $a_P: B_P \ra {\ub{\Z}}$ and
$g_V$. The map $\tilde{f}$ factorizes as $\mathcal Z_V \xra{f''} B'' \ra
B_P$ by the universal property of $D_V$. The main difference between $B''$ and 
$B_P$ lies in     their group of connected components.

\begin{proposition}\label{sette} The pair $(B'', f'')$ is the
Albanese scheme of $V$. \end{proposition}
 
\nd {\bf Proof.}    Consider a homomorphism $d': \mathcal Z_V \ra H'$
to a large group $H'$. 
By universality of $D_V$, the map $d'$ factorizes
uniquely as $jd$ where\ms 

(a)  $d: \mathcal Z_V \ra H$ is a homomorphism to a
large group $H$ inducing an isomorphism $D_V \iso
\pi_0(H)$; notice $H = H' \times_{\pi_0(H')} D_V$.
 
 and  
 
(b) $j: H \ra H'$ is a homomorphism of large groups.\ms

It suffices to show that the map $d$ factorizes via $(B'', f'')$. Let
$\mathcal H$ denote the sheaf of abelian groups on $S_{FL}$ corresponding
to the group scheme $H$. The composed map 
$d_1: \mathcal V \hra \mathcal Z_V \xra{d} \mathcal H$ of sheaves of sets
on $S_{FL}$ gives, by Yoneda's lemma, a $S$-morphism $d_1: V \ra H$
of schemes.  Since $V$ is connected, its image under $d_1$ is
contained in a connected component $Q$ of $H$. The neutral component
$C$  of $H$ is a semiabelian scheme. The scheme $Q$ has the structure
of a torsor under $C$ such that the pair of inclusions $Q \hra H$ and
$C\hra H$ combine into a morphism of torsors. In other words, there is a (injective)
morphism of group schemes $k:C_Q \ra H$.  Now we use
the fact that $f$ is universal for morphisms from $V$ to torsors under
semiabelian schemes. We get a
homomorphism 
$b: B\ra C$ and a morphism of torsors $h : P \ra Q$ such that $d_1=hf$. By 
Remark \ref{torso}, this is equivalent to specifying a homomorphism
$\tilde{h}: 
B_P \ra C_Q$ compatible with the augmentation to $\ub{\Z}$.  We have
shown that $d$ can be factorized as $k \tilde{h} \tilde{f}$. But, as we saw,
$\tilde{f}$ can be factorized via $f''$. This demonstrates the universal
property of the pair $(B'', f'')$. \qed\ms          

\begin{corollary}\label{old} Let $A_V$ be the Albanese scheme of a connected 
scheme $V$, as above. The neutral component of $A_V$ is the `` generalized 
Albanese variety"  of $V$ defined by Serre {\rm (\ref{serre})}. \end{corollary}

\begin{theorem}\label{alb2}  Let $X$ be a reduced algebraic scheme
over $S$. There exists a large group $A_X$ together with a map $u_X:
\mathcal Z_X \ra \mathcal A_X$ which is universal for homomorphisms
from $\mathcal Z_X$ to large groups. \end{theorem}

\nd {\bf Proof.} We have just seen that the Albanese scheme exists for connected
schemes. Now use Proposition \ref{enjoy}.\qed\ms

The identification of Corollary \ref{old} yields the 

\begin{corollary}

{\rm (i)} Given a torsor $P$ under a semiabelian scheme $B$, the
natural inclusion of  $P \xra{j} B_P$ induces a homomorphism $\tilde{j}:
\mathcal Z_P \ra B_P$.  The large group $B_P$, together with
the map $\tilde{j}$, is the Albanese scheme of $P$.

{\rm (ii)} If $X$ is proper, then $A_X^0$ is an abelian scheme.

{\rm (iii)} If $U$ is an open dense subscheme of $X$, then the induced map
$A_U \ra A_X$  is surjective.

 {\rm (iv) (homotopy invariance)} As usual, $\mathbb A^n$ denotes
 $n$-dimensional affine space over 
$S$.  Let 
$Y := X\times_S \mathbb A^n$. The natural projection $\pi: Y \ra X$ induces 
an isomorphism $\pi_*: A_Y \iso A_X$.

{\rm (v)}  For any semiabelian scheme $B$, one has $A_B = B \times_S
\ub{\Z}$.

{\rm (vi) (Kunneth)} Let $X:= Y \times_S V$ and let $\alpha: D_X \ra
D_Y \times_S D_V$ denote the natural map. One has $$ A_X = 
(A_Y \times_S A_V)\times _{(D_Y \times D_V)} D_X.$$ 

\end{corollary} 
\nd {\bf Proof.} We prove (iv), (v), (vi) and leave the rest to the
reader. For (iv), observe that $A_{\mathbb A^n} =~\ub{\Z}$:  there are no 
nonconstant morphisms from $\mathbb A^n$ to semiabelian
varieties. Statement (v) is a special case of (i) which follows from
Proposition \ref{sette}. As for (vi), it is clear that $\pi_0(A_X) =
D_X$. Now, use the fact that $A^0_X = A^0_Y \times_S A^0_V$ which
follows (using Corollary \ref{old}) from the corresponding property for the ``generalized Albanese
varieties'' of Serre (\ref{serre}) proved in \cite{Se2}.\qed\ms 

\begin{example}         A curve is a  scheme of pure 
dimension one over $S$.
Let $X$ be a normal integral curve.
 There is a unique
normal projective curve $X'$ 
containing $X$ as an open dense subscheme; let $R$ be the closed complement
of $X$ in $X'$. The semiabelian variety $A_X^0$ is the generalized
Jacobian of Rosenlicht \cite{Ro1} of $X'$ corresponding to the modulus $R$ 
(cf. \ref{rosen} and \cite{Se3} p. 4).\end{example}
 
\section{Special schemes} 

 Let $X$ be the open
complement of a divisor $E$ in  
a smooth and projective variety $X'$ of dimension $d$. Let $j: X \hra
X'$ denote the inclusion. Schemes such as $X$ are termed
\emph{special}. We fix these notations for the rest of this section.
We now turn to a discussion of  the ``generalized Albanese variety'' of
$X$ defined by Serre \cite{Se3}.

Let $w(E)$ denote the set of geometric points of $\pi_0(\tilde{E})$ i.e. $w(E)$
is the set of connected components of $\tilde{E}\times_S \bar{S}$ or,
what is the 
 same, the set of  irreducible components of $\bar{E}$. Given a set $W$, one
 defines the abelian group $\Z^W$ to be the group of integer valued
 functions on 
$W$ i.e. $\Z^W:= Maps(W,\Z)$. For a finite set $W$, the groups $\Z(W)$ and
 $\Z^W$ are naturally dual as abelian groups.
 
The \'etale group 
scheme $\ub{\Z}^E$ is defined by its geometric points: $\ub{\Z}^E(\bar{S}) = 
\Z^{w(E)}$. 
 If $Pic(\bar{X'}):= H^1(\bar{X'};\mathcal O_{\bar{X'}}^*)$ denotes
 the absolute Picard group of 
$\bar{X'}$, there is a natural {\bf g}-equivariant homomorphism
\begin{equation}
 \Z^{w(E)} \xra{b} Pic(\bar{X'}) \qquad d \mapsto \mathcal O(d);
\end{equation}  this provides 
a morphism\footnote{$Pic_{X'}$ is the Picard scheme of $X'$ (see \S \ref{picpro})}:
\begin{equation}\label{bee}
b: \ub{\Z}^{w(E)}  \ra D_{Pic^{red}_{X'}}  \end{equation}  
whose kernel is denoted $B$; set $I = B(\bar{S})$.  Under the
natural duality between $\Z^{w(E)}$ 
and $\Z(w(E))$, the group $I$ is dual to the subgroup $J$ of divisors of 
$X' \times_S \bar{S}$, supported on $E \times_S \bar{S}$, which are
algebraically   
equivalent to zero. Consequently, $I$ and $J$ have the same rank. The
{\it N\'eron-Severi} group 
$NS({\bar{X'}})$ of $\bar{X'}$ is the group of 
geometric points of $D_{Pic_{X'}}$. There is a natural morphism $B
\xra{v} Pic_{X'}^0$. 
 
 The abelian scheme $A^0_{X'}$ is the
``classical" Albanese  variety of $X'$. It is the maximal
abelian quotient of the semiabelian scheme $A^0_X$ (the ``generalized Albanese
variety" of Serre).

\subsection{Serre's Construction}\label{Serr}  Let us review the construction
of $A_X^0$,  via the   Picard scheme of $X'$,  
due to Serre \cite{Se3}.  
Denote by $v$ the  natural map from $B$ to the neutral 
component $Pic^{0}_{X'}$ of $Pic_{X'}$.
This defines the Picard 1-motive of $X$  

\begin{equation}\label{m1}
M^1(X):= [B \xra{v} Pic^{0,red}_{X'}]. \end{equation} 
The Cartier dual of this 1-motive is the Albanese 1-motive of $X$  
\begin{equation}\label{m2}
M_1(X):= [0 \ra G];\end{equation}  the group $G$ is  an extension
of $A_{X'}^0$  by the torus 
$T:= Hom(B,\gm)$.  

\begin{theorem}\label{serth} {\rm (J.-P. Serre)} The semiabelian
  scheme $G$ is  
naturally isomorphic to $A_X^0$ \cite{Se3}. \end{theorem}

\begin{corollary}\label{242} {\rm (a)}  The kernel of the natural 
 morphism $A_X^0 \ra A_{X'}^0$ is the torus $T= Hom(B, \gm)$. 
 
{\rm (b)} The dimension of the torus $T$ is the rank of 
the abelian group $ I$. 

{\rm (c)} the 1-motives (\ref{m1}) and (\ref{m2}) are 
independent of the compactification $X'$ of $X$. \end{corollary} 

Let $u : X' \to L'$ 
be a universal morphism  (see Remark \ref{serre}) where
$L'$ is a torsor for 
$A^0_{X'}$.  One has
the classical isomorphisms  
\begin{eqnarray}
Pic^{0}_{A^0_{X'}} \cong Pic^{0}_{L'} \cong
Pic^{0}_{X'},\\
P:= Pic^{0, red}_{X'} \qquad Pic^{0, red}_P \cong A^0_{X'}. \end{eqnarray}
 
Using this, one interprets the map $v$ as specifying certain
$\gm$-torsors on $A^0_{X'}$ (and on $L'$). The group scheme $G$ is the
geometric scheme corresponding to the total space of the $T$-bundle
(given by $v$) over $A^0_{X'}$. Likewise,  $L$ is the scheme over $L'$ corresponding
to the total space of the $T$-bundle. Using the universal map $u$, one
may pull this back to $X'$ to get a $T$-bundle over $X'$ with the
following 
property: given any element $b \in Hom(T,\gm)= B$, the
$\gm$-torsor over $X'$ obtained via $b$ corresponds
to the line bundle defined by $b$, viewed as a Cartier divisor on
$X'$. Observe the trivial fact that
 $X$ is disjoint from the support (contained in $E$) of all the elements of
$B(\bar{S})$ viewed as Cartier divisors on $X'$. Therefore, given
any $b$ as above, the $\gm$-torsor defined by $b$ has a canonical
trivialization over $X$: one has a rational section of the
$\gm$-torsor whose divisor is exactly $b$, and this is unique up to an
element of $H^0(X';\mathcal O^*_{X'})$. This
section can be used for the trivialization over $X$. Further, this can
be done in an 
uniform manner to derive a morphism $X \ra L$, which is unique up to
translations by elements of $T$. Serre \cite{Se3} shows that this
morphism is universal for maps to torsors under semiabelian schemes. This
extends, by universality of the sheaf $\mathcal Z_X$, to a morphism $\mathcal
Z_X \ra G_L$ where $G_L$ is the large group scheme defined by the
$G$-torsor $L$ (as in \ref{bps}). One deduces that this last morphism is the universal
morphism for $X$. 

The next theorem is true by definition of $M_1(X)$ and $M^1(X)$. But
the proof of the 
corresponding result  for simplicial smooth schemes (proved in \S
\ref{duldul}) is nontrivial. 

\begin{theorem}\label{duel}  The 1-motives $M_1(X)$ and $M^1(X)$  are  
dual. \end{theorem}

Let $Y$ be an open subscheme of $X$ whose complement in $X$ has
codimension greater than one. Let $B_Y$ be the group dual to the group
of divisors of $X'\times _S \bar{S}$ (supported on $(X' - Y)\times_S
\bar{S}$) which are algebraically 
equivalent to zero.  

\begin{corollary} The inclusion of $Y$ in $X$ induces an
isomorphism of the Albanese schemes $A_Y \iso A_X$.\end{corollary} 
  
\nd {\bf Proof.} Observe that $B$ and $B_Y$ are all the same.\qed\ms

\subsection{Hodge realization}  Just as a harbinger of later results,
we prove a result about the Hodge realization of $M_1(X)$. Consider an imbedding $\iota: k \hra
\C$ \footnote{The variety $X$ is always defined over a finitely
generated extension of $\Q$. So, we may assume that $k$ is finitely
generated over $\Q$.}. By base-change via $\iota$, we may and do assume (in this
subsection) that all our schemes
are over $\C$.   Since we are over an
algebraically closed field, there is no essential difference between a group
scheme and its torsors. Let  $u: X \ra A^0_X $ be the universal
morphism (uniquely determined by the choice of a base point). Recall
that $d$ is the dimension of $X$. 

Consider the long exact sequence in cohomology (with coefficients in
$\Z$) of the pair $(X',E)$:
$$ \ra H^{2d-2}(X')(d) \xra{\alpha} H^{2d-2}(E)(d) \ra
H^{2d-1}(X',E)(d) \ra H^{2d-1}(X')(d) \ra 0 {}$$
Poincar\'e duality gives $H^{2d-i}(X')(d) \cong H_i(X')$. By
Lefschetz duality, 
 we have that  $H^{2d-i}(X',E)(d) \cong H_i(X)$. The
image of $\alpha$ is the dual MHS of the kernel of $(H^0(E)(-1) \xra{\beta}
H^2(X'))$ where $\beta$ sends an element of $H^0(E)(-1)$ to its
Poincar\'e dual, an element of  $H^2(X')$.   The group $H^0(E)$ is
 the free abelian group $\Z^{w(E)}$ i.e. the group dual to $\Z(w(E))$.
  The map $b$ in (\ref{bee})  induces a map on the geometric points
  and there is a  
natural inclusion of   $D_{Pic_{X'}}(\C)$ in $H^2(X')$, via  the
exponential sequence. The composed map is the map $\beta$.
 The kernel of $\beta$
corresponds, under the duality of $\Z^{w(E)}$ and $\Z(w(E))$, to the 
subgroup $J$ of divisors of $X'$, supported on $E$, which are homologically
equivalent to zero\footnote{For divisors, homological
equivalence coincides with algebraic equivalence.}.  Hence we find
that Ker$\beta = I$.  If $I^* = Hom_{\mathbb Z}(I, \mathbb Z)$, 
then Im~$\alpha \cong  Hom_{\mathbb Z}(I(-1), \mathbb Z) = I^*(1)$; both
are MHS of type ($-1, -1$). So we obtain the exact sequence  $$ 0 \ra I^*(1) \ra H_1(X)
\ra H_1(X') \ra 0.{}  $$ 
\begin{theorem} {\rm (Hodge realization)}  

{\rm (i) } The map 
$$u_*: H_1(X; \Z)/{tor} \ra H_1(A_X^0; \Z) =: \mathfrak T(M_1(X)){} $$
of mixed Hodge structures is an isomorphism.

{\rm (ii) } There is a natural  isomorphism 
$$          \mathfrak T (M^1(X)) \xra{u^*} H^1(X;\Z(1))   {} $$ 
of mixed Hodge structures. \end{theorem}

\nd {\bf Proof.} Recall from (\cite {Se1} Chap. VI, \S 14, Prop. 13,
p.128) that $u_*$ is surjective; universal morphisms are maximal in
the sense of {\it loc. cit}.   

Consider the following commutative diagram:
$$
\begin{CD}
 0  @>>> I^*(1) @>>> H_1(X)/{tor}     @>>> H_1(X')/{tor}      @>>> 0  \\
 @.      @VVV        @V{u_*}VV         @V{u'_*}VV        @.  \\
 0  @>>> H_1(T) @>>> H_1(A_X^0) @>>> H_1(A_{X'}^0) @>>> 0   \\
\end{CD} {} 
$$
All the vertical maps are maps of MHS because they are induced by morphisms. 
It is a classical result (see p. 553 of \cite{gh}) that $u'_* : H_1(X')/{tor} \ra H_1(A_{X'}^0)$ is an
isomorphism.  By the snake lemma, $I^*(1) \ra H_1(T)$ is surjective. But  the 
rank of $H_1(T)$ is the same as the dimension 
of $T$ which equals the rank of $I$ , as remarked earlier.  Since they both 
have the same rank, the left vertical map is an isomorphism. Thus, we
obtain that 
the map $u_* : H_1(X)/{tor} \ra H_1(A_X^0)$ is a isomorphism of MHS. 
By definition,  $\mathfrak T(M_1(X))$ is the mixed Hodge
structure $H_1(A^0_X; \Z)$  \cite{h}, 10.1.3 p. 54.
We turn to (ii).
The universal coefficient theorem says that the groups $H_1(X;\Z)/{tor}$ and
$H^1(X;\Z)$ are naturally dual. So $H^1(X)(1)$ and $H_1(X)/{tor}$ are dual as
 elements of $\mh$, i.e., $H^1(X)(1) = Hom_{MHS}(H_1(X),
\Z(1))$. Since duality of  1-motives is compatible with duality
in $\mh$ (\cite{h}, \S 10.2), the result
follows from the first part.\qed\ms

\begin{remark} One cannot expect the previous result to be true for an 
arbitrary complex algebraic  variety (thereby justifying the entrance
  of 1-motives in this context): the weights in the $H_1$ of a 
semiabelian variety are -2 and -1, whereas the weights that can and do occur in
the $H_1$ of an arbitrary variety are -2, -1 and 0 (\cite{hp} 7.1,
  8.3). A noncompact  nonrational curve with  
nodes provides an easy example of the last property.\end{remark}

\section{Picard schemes}\label{picpro}

In this section, we extend the classical construction of the Picard
scheme to a simplicial setting. The arguments are based on spectral
sequences. 

  Let $f: V \ra S$ be a proper reduced scheme of finite type over $S$.

We define $\mathcal R^1f_*\mathcal O^*_V$ on $S_{FL}$ to be the sheaf
associated with the presheaf $T \mapsto H^1(V \times_S T; \mathcal O^*_{V
\times_S T})$. By a theorem of J. Murre and F. Oort (\cite{BLR}, \S 8.2,
Theorem 3, p. 211), the sheaf $\mathcal R^1f_*\mathcal O^*_V$ is representable
by a locally algebraic group scheme $Pic_V$. One refers to
$Pic_V$ as the Picard scheme of $V$. One obtains that 
$$H^1(\bar{V}; \mathcal O^*_{\bar{V}}) = Pic_V(\bar{S}).$$ 

We recall the well known \footnote{The crystalline analog may be found
in \cite{ill} II 3.11 and  \cite{ill2} pp.46-47.} 
\begin{proposition}\label{derpic} 
{\rm (i) (de Rham realization)} Let $k$ be of characteristic
zero. There is a canonical isomorphism {\rm (\cite{memu}, \cite{Me},
\cite{h} \S 10.1)}   
$$ \mathfrak T_{dR}([0 ~\ra Pic^0_V]) \cong H^1_{dR} (V). {} $$

{\rm (ii)(Hodge realization)} Let $k = \C$. The exponential sequence
and G.A.G.A. furnish a natural isomorphism
$$H^1(V)(1) \iso \mathfrak T([0~\ra Pic^0_V]). {} $$
\end{proposition}

\subsection{The Picard scheme of proper simplicial schemes.}
 Let $Z_{\bu}$ be a simplicial scheme. The {\it absolute} Picard
group $Pic(Z_{\bu})$ of $Z_{\bu}$ is the group $H^1_{Zar}(Z_{\bu};\mathcal
O^*_{Z_{\bu}})$ where $\mathcal O^*_{Z_{\bu}}$ is the sheaf of units of
the structure sheaf $\mathcal O_{Z_{\bu}}$ of $Z_{\bu}$. We often use
the notation $\mathcal O^*$ to indicate the sheaf of units of a simplicial
scheme (specified by the context). 

By the definition of cohomology in the
simplicial setting (\cite{h}, pp. 12-14, esp. 5.2.3 and 5.2.7,
\cite{gi}, Ex. 1.1, p.7),
elements of the group $Pic(Z_{\bu})$ correspond to isomorphism classes
of pairs $(\mathcal L, \alpha)$ where    

(1) $\mathcal L$ is an invertible sheaf on $Z_0$ such that $\delta^*_0
\mathcal L$ is isomorphic to $\delta^*_1\mathcal L$;    

(2) $\alpha$ is an isomorphism $ \delta_0^* {\mathcal L} \iso \delta_1^*
{\mathcal L}$ on  $Z_1$ satisfying a cocycle condition on $Z_2 $, i.e.,
the following diagram is commutative:       

$$
\begin{CD}
\delta_2^* \delta_0^* \mathcal L @>{\delta_2^* \alpha}>> \delta_2^* \delta_1^*
\mathcal L @ =
\delta_0^* \delta_0^* \mathcal L \\
@|    @.   @V{\delta_0^* \alpha}VV \\
\delta_1^* \delta_0^* \mathcal L @>{\delta_1^* \alpha}>> \delta_1^*
\delta_1^* \mathcal L @=
\delta_0^* \delta_1^* \mathcal L \\ 
\end{CD} {} 
$$

The maps $\delta_0, \delta_1 : Z_1 \ra Z_0$ and $\delta_0, \delta_1,
\delta_2: Z_2 \ra Z_1$ form part 
of the simplicial structure of the scheme $Z_{\bu}$ (\cite{h} III,
5.3.6, pp. 15-16).  We write the cocycle condition  as
\begin{equation}\label{rowss}
\delta_2^*(\alpha)\delta_0^*(\alpha) = \delta_1^*(\alpha).{} 
\end{equation}

Let $V_{\bu}$
be a {\it proper} reduced simplicial scheme. The structure morphism of
$V_n \ra S$ is denoted  by $f_n$ and the structure map $V_{\bu} \ra
S$ is denoted by $f$.       
For any scheme $T$ over $S$, we get a simplicial scheme
$V_{\bu}\times_S T$. By \cite{h} III, 5.2.3.2, p. 13, there is a
spectral sequence 
\begin{equation}\label{specs}
E_{1,T}^{p,q} = H^q(V_p \times_S T;\mathcal O^*_{V_p
\times_S T}) \Longrightarrow H^{p+q}(V_{\bu} \times_S T;\mathcal O^*).
\end{equation}
We define $\mathcal R^1f_*\mathcal O^*$ to be the sheaf on
$S_{FL}$ associated with the presheaf $E^1$:
 $T \mapsto H^1(V_{\bu} \times_S T; \mathcal O^*)$. 
The first result that one has about this sheaf is the 
{}
\begin{theorem}\label{repic} The sheaf $\mathcal R^1f_*\mathcal O^*$ on
$S_{FL}$ is representable by a locally algebraic group scheme (denoted
$Pic_{V_{\bu}}$). \end{theorem}

\nd {\bf Proof.} We define $\mathcal R^if_{n*}\mathcal O^*_{n}$ to be the
sheaf  on $S_{FL}$ 
associated with the presheaf: 
$$T \mapsto H^i(V_n \times_S T; \mathcal O^*_{V_n \times_S T}) {} $$  
The sheaves $\mathcal R^0f_{n*}\mathcal O^*_{n}$ are representable by tori.
By the classical theorem of Murre and Oort (\cite{BLR}, \S 8.2,
Theorem 3, p. 211), the sheaves $\mathcal
R^1f_{n*} \mathcal O^*_{n}$  are representable by  locally algebraic group
schemes. 

The spectral sequence (\ref{specs})  gives rise to an exact sequence
(natural in $T$):
\begin{equation}\label{lowly}
0 \ra E^{1,0}_{2,T} \ra H^1(V_{\bu} \times_S T;\mathcal O^*) \ra
E^{0,1}_{2,T} \ra E^{2,0}_{2,T}.\end{equation}

Let $E^{p,q}$ be the presheaf defined by $T \mapsto
E^{p,q}_{2,T}$ and let $\mathcal E^{p,q}$ 
be the sheaf 
associated with the presheaf $E^{p,q}$.

One has an exact sequence of sheaves on $S_{FL}$:
$$ 0 \ra \mathcal E^{1,0} \ra  \mathcal R^1f_{*}\mathcal O^*  \ra
\mathcal E^{0,1} \xra{d_2} \mathcal E^{2,0}.{} $$

\nd  The sheaf $\mathcal E^{1,0}$ is the homology of
the complex of sheaves: 
$$ \mathcal R^0f_{0*}\mathcal O^*_{0}  \xra{\delta^*_0 - \delta^*_1} \mathcal
R^0f_{1*}\mathcal O^*_{1} \xra{\delta^*_0 -\delta^*_1 + \delta^*_2} \mathcal
R^0f_{2*}\mathcal O^*_{2}. {} $$ 
Each of the sheaves in the above complex is representable by a torus.
Therefore the sheaf $\mathcal E^{1,0}$ is representable by an affine group
scheme $C'$ whose neutral component $C$ is
a torus.   

- The same argument demonstrates the representability
of the sheaf $\mathcal E^{2,0}$ by an affine group scheme $W$ with
$W^0$ a torus.
 
- The sheaf $\mathcal E^{0,1}$ is the sheaf $Ker~(\mathcal
R^1f_{0*} \mathcal O^*_{0} \xra{\delta^*_0 - \delta^*_1} \mathcal R^1f_{1*}
\mathcal O^*_{1})$. Since each of these last sheaves is representable by
locally algebraic group schemes, we find that the sheaf $\mathcal E^{0,1}$
is also representable by a locally algebraic group scheme, denoted $K$. 

- The homomorphism $d_2$  corresponds to a morphism of
locally algebraic group schemes
$K \xra{d} W. {} $  The kernel of $d$  is a
locally algebraic group scheme $Q'$ representing the sheaf
$Ker(d_2)$. 
Thus, in  the exact sequence 
$$ 0 \ra \mathcal E^{1,0} \ra  \mathcal R^1f_{*}\mathcal O^*  \ra
Ker(d_2) \ra 0,$$
the left term is representable by an affine scheme, and the right term
is representable. The representability of the sheaf $ \mathcal 
R^1f_{*}\mathcal O^* $ follows from the next proposition. \qed\ms 

\begin{proposition}\label{repo}
Let $0 \rightarrow F \rightarrow G \rightarrow H \rightarrow
$ be an exact sequence of sheaves of abelian groups on
$S_{FL}$. Suppose $F,H$ are representable, and that $F \rightarrow
S$ is an affine morphism (i.e. the scheme representing $F$ is affine
over $S$). Then $G$ is representable (necessarily by a commutative
group scheme).
\end{proposition}

\nd {\bf Proof.} This is Proposition 17.4 on page III.17-6 of
\cite{oort}. \qed\ms

The Picard scheme $Pic_{V_{\bu}}$ of $V_{\bu}$ is an extension of $Q'$ by the
(reduced) affine group scheme $C'$. We now recall the (\cite{DG} Corollary
2.2, II \S 5, no. 2, p. 288):  

\begin{theorem}\label{reddy} Let $G$ be a locally algebraic group scheme over a field
$F$. If the field $F$ is perfect, then the 
reduced scheme $G^{red}$ is a smooth group scheme.\end{theorem}

This is the main reason why we assume that $k$ is perfect. Other
reasons include 

\nd (i) the structure theorems for commutative group schemes are
particularly nice over a perfect field (Theorems 1
and 2 of \S 9.2 in \cite{BLR}).  

\nd (ii) Over a perfect field, a curve $X$ is smooth if and only if it
is normal. These two notions are not equivalent over arbitrary fields
(related to the distinction between smoothness and regularity of
schemes over imperfect fields; cf. Exercise 10.1 of \cite{Ha}). 

\nd (iii) Grothendieck has expressed doubts about the suitability of
1-motives over imperfect fields (see the appendix to
\cite{jann}). This is based on the fact that the 1-motives associated
with a (singular) curve need not be compatible with base change.   

\nd (iv) our eventual interest is in constructing 1-motives for general
varieties by using simplicial smooth varieties. Our construction (\S
\ref{positive}) is based on the results
of \cite{dj} which are valid only over perfect fields (page 1 of
\cite{dj}).\medskip 
  
\begin{corollary}\label{pic} Let $V_{\bu}$ be a proper reduced
simplicial scheme. Assume that $V_0 \times_S \bar{S}$ is normal i.e.
$V_0$ is geometrically normal. The  group scheme
$Pic^{0,red}_{ V_{\bu}}$ is a semiabelian scheme.
\end{corollary}

\nd {\bf Proof.} The sheaf $\mathcal R^1f_{0*} \mathcal O^*_{0}$  is
representable by a group scheme whose reduced subscheme is an
      {\it abelian} locally algebraic  scheme \cite{G1}, Theorem 2.1,
p. 11. In other words, the neutral component of $K^{red}$  is an
abelian scheme.  Recall 
that every morphism from an 
abelian variety to an affine group scheme is constant; in particular,
a homomorphism must necessarily be the zero map. Therefore, the neutral 
component of $K^{red}$
must be in the kernel of the morphism $d$. In other words, if $Q$
denotes the neutral component of $(Q')^{red}$, the
inclusion of $Q$ in $K^{0,red}$ is an equality. 
We obtain that
$Pic^{0, red}_{V_{\bu}}$ is a semiabelian scheme.\qed\ms

In particular, Corollary \ref{pic} applies  to a {\it smooth} proper
simplicial scheme.  

\begin{proposition}\label{pic2} Let $V_{\bu}$ be a proper reduced simplicial
scheme. One has the following natural isomorphisms:

(0) $ H^1(V_{\bu}\times_S \bar{S}; \mathcal O^*) \iso Pic_{V_{\bu}}(\bar{S})$

(1) $ H^1(V_{\bu}; \mathcal O ) \iso Lie~Pic^0_{V_{\bu}} = 
Lie~Pic_{V_{\bu}}$

 (2) ($k = \C$)  $ H^1(V_{\bu}; \mathcal O^*) \iso  H^1(V^{an}_{\bu};
 \mathcal O^*),  
 \qquad H^1(V_{\bu}; \mathcal O) \iso  H^1(V^{an}_{\bu}; \mathcal O)
 $

 (3) ($k= \C$)  $ H^1(V_{\bu}; \Z(1)) \iso  H_1(Pic^0_{V_{\bu}};
 \Z)$. 
\end{proposition}

\nd {\bf Proof.} The  isomorphisms in (2) is a natural consequence of
 G.A.G.A.  applied to $E^{p,q}_{2,S}$ of (\ref{lowly}).  
Assertion (0) is clear. The proof 
of (1) follows the proof of Theorem 1, \S 8.4 on p. 231 in \cite{BLR}. 
In the exact sequence (induced by the exponential sequence on $V_{\bu}$):
$$  H^1(V_{\bu}; \Z(1)) \xra{\partial}  H^1(V^{an}_{\bu}; \mathcal O)
\xra{exp}  H^1(V^{an}_{\bu}; \mathcal O^*) \ra  H^2(V_{\bu};\Z(1)),
{} $$ 
the map $\partial$ is injective since 
$exp:\C \ra \C^*$ is surjective. 
By parts (1) and (2), this is enough to prove (3). \qed\ms

\section{The Albanese and Picard 1-motives}\label{ilbee}
Let us fix a simplicial scheme $X_{\bu}$, the 
complement of a divisor $E_{\bu}$ with normal crossings in a smooth
projective simplicial scheme $X'_{\bu}$. In this section, we define
the Albanese and Picard 1-motives of $X_{\bu}$.  We
prove a
conjecture of Deligne that the Hodge realization of the Picard
1-motive of $X_{\bu}$ (over $\C$) is $H^1(X_{\bu}; \Z(1))$. 

\subsection{The Picard 1-motive.}\label{alss}

Let $R_n$ be the large group scheme  $\Z^{E_n}$. Put $R"_n:= Ker(b_n: R_n
\ra D_{Pic_{X_n}})$.    We define $$R'
:=~Ker(R_0 \xra{\delta_0^* - \delta_1^*} R_1).{} $$ 

Any element $r$ of $R'(\bar{S})$
determines an invertible sheaf $\mathcal O(r)$ on $X'_0 \times_S \bar{S}$.
There is a
canonical isomorphism $\alpha_r: \delta^*_0 \mathcal O(r) \iso
\delta^*_1 \mathcal O(r)$ on $X'_1 \times_S \bar{S}$; the isomorphism
$\alpha_r$ satisfies the  
cocycle
condition $$p_2^*(\alpha_r)p_0^*(\alpha_r) = p_1^*(\alpha_r). {} $$
We have a homomorphism $$R'(\bar{S}) \xra{h} Pic_{X'_{\bu}}(\bar{S})=
Pic(X_{\bu} \times_S \bar{S}) \qquad h'(m)=
([\mathcal O(r)], \alpha_r). $$  Let $R' \xra{h'} Pic_{X'_{\bu}} $ be the corresponding
homomorphism of group schemes. Composing with the natural morphism 
$Pic_{X'_{\bu}} \ra NS({X'_{\bu}}):= \pi_0(Pic_{X'_{\bu}})$, one gets
a homomorphism  
$$R' \xra{h"} NS({X'_{\bu}}) {} $$ whose kernel is denoted
$R$. Another description of $R$ is 
$$R = Ker(R"_0 \xra{\delta_0^* - \delta_1^*} R"_1).{} $$ Recall
from Corollary \ref{pic} that  
$P:= Pic^{0, red}_{X'_{\bu}}$ is a semiabelian scheme.

\begin{definition} The motivic $H^1$ (or the Picard 1-motive) of $X_{\bu}$ is: $$M^1(X_{\bu})~:= ~[R \xra{h} P]. {} $$ 
\end{definition}

\subsection{The Albanese 1-motive} 
We put $A_i= A_{X_i}$ and $D_i= D_{X_i}$; note that $D_i = D_{A_i}$. 
By virtue of the simplicial structure of $X_{\bu}$, the groups $A_i$,
$A^0_i$, $D_i$  assemble to form simplicial groups $A_{\bu}$,
$A^0_{\bu}$ and $D_{\bu}$; they  fit into an exact  sequence of
simplicial group schemes 
$$ 0 \ra A^0_{\bu}\ra A_{\bu} \ra D_{\bu} \ra 0. $$

There is a natural functor (usually termed normalization) \cite{may}
which transforms a  simplicial commutative group into a  homological
complex; the simplicial 
maps $\delta_i$ are replaced by the differentials $d:= \Sigma_i (-)^i
\delta_i$. Thus, the previous sequence  can be
interpreted as an exact sequence of complexes of locally algebraic
group schemes 
(concentrated in nonnegative degrees); as such, one gets a long exact
sequence of homology.  We set 
$$G:= H_0(A^0_{\bu}) = \frac{A^0_0}{d(A^0_1)}, \qquad N:=
H_1(D_{\bu}) = \frac{Ker(d:D_1 \to D_0)}{Im(d: D_2 \to
D_1)}.$$ One has the boundary map $\partial: N \to G$.
 
 \begin{definition} The Albanese 1-motive (or the motivic $H_1$) of $X_{\bu}$ is  $$M_1(X_{\bu})
 := [N\xra{\partial} G].$$\end{definition} 

Strictly speaking, $N$ has to be torsion-free for $M_1(X_{\bu})$ to be
a 1-motive; but we can still make sense of the dual 1-motive of
$M_1(X_{\bu})$. So the torsion in $N$ does not affect matters much. 

\begin{remark}\label{constance} (Constant schemes) If $X_{\bu}$ happens to be a constant
simplicial scheme, i.e., if $X_i = X_0$ for all $i > 0$ and each of the 
simplicial maps is the identity, then $M_1(X_{\bu}) = M_1(X_0)$ where
$M_1(X_0)$ of the special scheme $X_0$ is defined as in
(\ref{m1}). \end{remark} 

\begin{theorem}\label{picho} {\rm (Deligne's conjecture)}\footnote{This has been proved (independently) by
Barbieri-Viale and Srinivas \cite{bs2}; it was announced in
\cite{bs}.}  Let $k = \C$. 
One has a natural
isomorphism $$ H^1(X_{\bu};\Z(1)) \ra \mathfrak T (M^1(X_{\bu})).{} $$ 
\end{theorem}
  
\nd {\bf Proof.}  As before, we set $P:=
Pic^0_{X'_{\bu}}$; it is a reduced scheme. 

The main step is the construction of a homomorphism $\beta$ from
$H^1(X_{\bu};\Z(1))$ to Lie~$P$, the Lie algebra of $P$; note that
Lie~$P$ is $H^1(X'_{\bu}; \mathcal O_{X'_{\bu}})$. For this, we closely
follow Deligne (\cite{h}, 10.3.9-10.3.15, p. 71-74); the modifications
are minor. 

\begin{lemma}\label{subs} {\rm (i)} Let $j^m_*\mathcal O^*_X$ be the subsheaf of
  meromorphic functions of 
$j_*\mathcal O^*_X$ on $X'_{\bu}$. One has
$$H^1(X_{\bu};\Z(1)) = \mathbb H^1(X'_{\bu};[\mathcal O_{X'} \xra{{exp}}
j^m_*\mathcal O^*_{X}]). {} $$
{\rm (ii)} $$H^1(X_{\bu};\C) = \mathbb H^1(X'_{\bu};[\mathcal O_{X'} \xra{d}
\Omega^1_{X'_{\bu}}(log~ E)]). {} $$
{\rm (iii)} The inclusion of $H^1(X_{\bu};\Z(1))$ in $H^1(X_{\bu};\C)$ is
defined by the morphism of complexes
$$ \begin{CD}
\mathcal O_{X'} @>{exp}>> j^m_*\mathcal O^*_X \\
@| @V{df/f}VV\\
\mathcal O_{X'} @>d>> \Omega^1_{X'_{\bu}}(log~ E).\\
\end{CD}
{} 
$$

\end{lemma}

\nd {\bf Proof.} Consider the commutative diagrams
$$
\begin{CD}
H^1(X_{\bu};\Z(1)) @>{\sim}>> \mathbb H^1(X_{\bu};[\mathcal O_X @>{exp}>>
\mathcal O^*_X])  @>{\sim}>>
\mathbb H^1(X'_{\bu};[j_*\mathcal O_X @>{exp}>> j_*\mathcal O^*_X])\\
@VVV  @VVV  @. @VVV\\
H^1(X_{\bu};\C) @>>> \mathbb H^1(X_{\bu};[\mathcal O_{X} @>d>>
\Omega^1_{X_{\bu}}]) @<<<
\mathbb H^1(X'_{\bu};[j_*\mathcal O_X @>>> j_*\Omega^1_{X_{\bu}}])\\
\end{CD}
{} 
$$

$$
\begin{CD}
\mathbb H^1(X'_{\bu};[j_*\mathcal O_X @>{exp}>> j_*\mathcal O^*_X])
@<\sim<< \mathbb H^1(X'_{\bu};[\mathcal O_{X'} @>{exp}>> j^m_*\mathcal O^*_X])\\
 @VVV @. @VVV @. \\
\mathbb H^1(X'_{\bu};[j_*\mathcal O_X @>>> j_*\Omega^1_{X_{\bu}}]) @<{\sim}<< \mathbb
H^1(X'_{\bu};[\mathcal O_{X'} @>>> \Omega^1_{X'_{\bu}}(log~ E)]).\\
\end{CD}
{} 
$$

The horizontal maps in the first diagram are isomorphisms due to

(1) the exact sequence of complexes on $X'_{\bu}$:
$$ 0 \ra \Z(1) \ra \mathcal O_{X_{\bu}} \xra{{exp}} \mathcal O^*_{X_{\bu}}
\ra 0. {} $$
(2) the fact that the inclusion of complexes
$$ 0 \ra \C \ra [\mathcal O_{X} \xra{d} \Omega^1_{X_{\bu}}] {} $$
is a quasi-isomorphism in degree less than 2. More generally, the
inclusion of the truncated De Rham complex $\Omega^{<n}_{V_{\bu}}$ in
the De Rham complex of a smooth simplicial scheme $V_{\bu}$ induces a
quasi-isomorphism in degree less than $n$.

(3) the fact that  $\mathcal R^1j_*\mathcal O_{X_{\bu}} =0$ (\cite{h}
    3.1.7). 

The horizontal maps in the second diagram are
isomorphisms because  the morphism of complexes which define them are
quasi-isomorphisms.\qed\ms

 The group $\mathbb H^1(X'_{\bu}; [\mathcal O_{X'_{\bu}}^* \ra \Omega_{X'_{\bu}}^1 (log~ E)])$ can be interpreted
as  the group of isomorphism classes of triples $(\mathcal L, \alpha_{\mathcal
L}, \omega_{\mathcal L})$, with $\mathcal L$ an invertible sheaf on
$X'_{\bu}$ (thereby determining $\mathcal L_0$ on $X'_0$), $\alpha_{\mathcal L}$ is an
isomorphism $ \delta ^*_0(\mathcal L_0) \iso \delta^*_1(\mathcal L_0)$, and 
$\omega_{\mathcal L}$ is an {\it integrable} connection on $\mathcal
L$, regular on $X_{\bu}$ with 
at most simple poles along the divisor $E_{\bu}$ (cf. \cite{Me} 2.5,
p. 364 or \cite{RSS}, 1.5, p. 47),
whereby, using the \emph{connections on invertible sheaves = one-forms}
dictionary, 
one can identify $\omega_{\mathcal L}$ with an element of
$H^0(X'_{\bu};\Omega^1_{X'_{\bu}}(log~ E))$. 

If $a$ denotes the projection $X'_{\bu}$ to $S$, then the $S_{FL}$ sheaf 
$R^1a_*[\mathcal O_{X'_{\bu}}^* \ra \Omega_{X'_{\bu}}^1 (log~ E)]$ is
representable by a group scheme,
 denoted 
$M^{\natural}_{X_{\bu}}$, whose group of $\C$-points is the group 
$\mathbb H^1(X'_{\bu}; [\mathcal O_{X'_{\bu}}^* \ra \Omega_{X'_{\bu}}^1 (log~
E)])$. The group scheme 
$M^{\natural}_{X_{\bu}}$ is an 
extension of a subgroup (containing $P$) of $Pic_{X'_{\bu}}$  by the
additive group scheme  
representing the $S_{fppf}$ sheaf $R^0a_*\Omega^1_{X'_{\bu}}(log~ E)$. The
representability of $R^1a_*[\mathcal O_{X'_{\bu}}^* \ra
\Omega_{X'_{\bu}}^1 (log~ E)]$ follows by applying Proposition
\ref{repo} to the exact sequence
$$ 0 \to  R^0a_*\Omega^1_{X'_{\bu}}(log~ E) \to R^1a_*[\mathcal O_{X'_{\bu}}^* \ra
\Omega_{X'_{\bu}}^1 (log~ E)] \to  R^1a_*\mathcal O_{X'_{\bu}}^*;$$
details may be found in \cite{me2}. 

Consider an element $r$ of $R(\C)$, i.e. a divisor $r_0$ supported on $E_0$,
algebraically equivalent to zero, such that $\delta^*_0(r_0) =
\delta^*_1(r_0)$ (equality of divisors on $X'_1$). The invertible sheaf
$\mathcal O(r_0)$ on $X'_0$ extends to an invertible sheaf $\mathcal O(r)$ on $X'_{\bu}$,
equipped with a canonical isomorphism $\alpha_r: \delta^*_0(\mathcal O(r_0))
\cong \delta^*_1(\mathcal O(r_0))$, and an integrable connection
$\omega_r$ on $\mathcal O(r)$, regular on $X_{\bu}$ with atmost simple
poles along the divisor $E_{\bu}$. The connection $\omega_{r}$ or,
equivalently the corresponding one-form, may be described as follows:
If $f_i$ are defining equations for the divisor $r_0$ on any open
covering {$U_i$} of $X'_0$, then $\omega_r= d log
(f_i) = d(f_i)/{f_i}$ on $U_i$.      
One defines a map
\begin{equation}\label{natmap}
R \ra M^{\natural}_{X_{\bu}} \qquad r \mapsto (\mathcal O(r), \alpha_r
, \omega_r)
\end{equation}

The group $H^1(X_{\bu};\C) =\mathbb H^1(X'_{\bu}; [\mathcal O_{X'_{\bu}} \ra
\Omega_{X'_{\bu}} (log~ E)])$ is the Lie algebra of
$M^{\natural}_{X_{\bu}}$. It can be identified with the isomorphism
classes of triples $(L, \alpha_L, \omega_L)$, for $L$
 an $\mathcal O_{X'_{\bu}}$-torsor (thereby determining a $\mathcal
O_{X'_0}$-torsor $L_0$), an isomorphism $\alpha_L: \delta^*_0(L_0) \cong \delta^*_1(L_0)$ and
$\omega_L$ an integrable connection, as above, on $L$.

The group $H^1(X_{\bu};\Z(1)) =\mathbb H^1(X'_{\bu}; [\mathcal O_{X'_{\bu}} \xra{{exp}}
j^m_*\mathcal O^*_{X_{\bu}}])$ is the group of isomorphism classes of
triples $(L, \alpha_L, \beta_L)$, with $L$ an $\mathcal
O_{X'_{\bu}}$-torsor (determining an $\mathcal O_{X'_0}$-torsor $L_0$), an
isomorphism $\alpha_L: \delta^*_0(L_0) \iso \delta^*_1(L_0)$, and
$\beta_L$ an isomorphism of the invertible sheaf $exp(L)$ with an
invertible sheaf $\mathcal O(r)$ (here $r \in R(\C)$). The map
$Aut(L) = \C \ra Aut(exp(L)) = \C^*$ is surjective. 
This allows one to identify $H^1(X_{\bu};\Z(1))$ with the group of
triples $(t, \alpha_t, r)$ for $t$ an isomorphism class of $\mathcal
O_{X'_{\bu}}$-torsor (determining a $\mathcal O_{X'_0}$-torsor $t_0$) and
an isomorphism $\alpha_t:\delta^*_0(t_0) \cong \delta^*_1(t_0)$,
i.e. an element of $Lie~ Pic_{X'_{\bu}}$, and $r\in R(\C)$, defining a
divisor concentrated on $E_{\bu}$ with $exp(t)$ as image in $P=
Pic^0_{X'_{\bu}}(\C)$.
In other words, we have defined an isomorphism (of abelian
groups)\footnote{See the definition of $\mathfrak T$ \cite{h} \S
10.1.3.1.}  
\begin{equation}\label{isoso}
H^1(X_{\bu};\Z(1)) \iso \mathfrak T(M^1(X_{\bu})). 
\end{equation}
The following lemma tells us that this is compatible with the weight filtration and
the Hodge filtration. 
{}
\begin{lemma}\label{seedy} {\rm (i)} One has
\begin{eqnarray*}
W_1(H^1(X_{\bu};\Z)) = Im((H^1(X'_{\bu};\Z)) \xra{{j^*}}
(H^1(X_{\bu};\Z))),\\
W_0((H^1(X_{\bu};\Z))) = Im(Ker(H^1(X'_{\bu};\Z) \ra
H^1(X_0;\Z)) \xra{{j^*}} H^1(X_{\bu};\Z)).
\end{eqnarray*} 
{\rm (ii)} The spectral sequence defined by the ``stupid'' filtration
of $[\mathcal O_{X'_{\bu}} \ra\Omega^1_{X'_{\bu}} (log~ E)]$ degenerates
and converges to the Hodge filtration on $H^1(X_{\bu};\C)$. \end{lemma}
\nd {\bf Proof.} 
For the first, let us check that $H^1(X_{\bu}~{\rm  mod}~ X'_{\bu}
;\Z)$ is pure of weight two. This follows from  the exact
sequence 
$$ 0 \ra  H^1(X'_{\bu}; \Z(1)) \xra{{j^*}}   H^1(X_{\bu}; \Z(1)) \xra{{Div}} 
H^0(E_{\bu}; \Z) \xra{d} H^2(X'_{\bu};\Z(1)) {} $$
(with a step of the Tate twist $\Z(-1)$, since $H^0(E_{\bu};\Z)$ is
pure of weight zero), part of the cohomology sequence of the exact
sequence of complexes on $X'_{\bu}$  
$$ 0 \to [\mathcal O_{X'} \xra{exp}
\mathcal O^*_{X'}] \to [\mathcal O_{X'} \xra{exp} j^m_*\mathcal O^*_X] \to
[0 \to \Z(E_{\bu})] \to 0;$$ 
the last map is the divisor map $j^m_*\mathcal O^*_X \to \Z(E_{\bu})$. 

 The second assertion in (i) follows from examining 8.1.19.1 and
8.1.20 of \cite{h}, p.35. 
As for (ii), it is clear from the definition 3.1.11, 8.1.8, 8.1.12 of
\cite{h} of the Hodge filtration. \qed\ms 

Back to the proof of Theorem \ref{picho}:  Combining the commutativity of the
diagram
$$
\begin{CD}
 H^1(X_{\bu}; \Z(1)) @>{\sim}>> \mathfrak T_{\Z}(M^1(X_{\bu})) @>>> Lie~ P =
H^1(X'_{\bu};\mathcal O_{X'_{\bu}})\\
@VVV @. \qquad @|\\
H^1(X'_{\bu}; \C) @>{\sim}>> Lie~ M^{\natural}(X_{\bu}) @>>> Lie~ P =
H^1(X'_{\bu};\mathcal O_{X'_{\bu}})\\
\end{CD}
{} 
$$ with part (ii) of Lemma \ref{seedy}, we
see that the map in (\ref{isoso}) is compatible with the Hodge
filtration. \qed\ms

\section{A Duality theorem}\label{duldul}
We retain the notations of \S
\ref{ilbee}.
In this section, we prove the duality of the Albanese and the Picard
1-motives (conjectured by Lichtenbaum). Recall the exact sequence
(here $P:= Pic^{0, red}_{X'_{\bu}}$) 
$$ 0 \to C' \to P \xra{\pi} Q \to 0.$$ We set $a: X'_{\bu} \ra S$ to be
the natural projection; we use $a_i: 
X'_i \ra S$  for each of the components. 

\subsection{The duality theorem.}

\begin{theorem}\label{duality} {\rm (Lichtenbaum's conjecture)} 
The 1-motives $$M_1(X_{\bu}):= [N \xra{\partial} G], \quad
M^1(X_{\bu}):= [~R \xra{h} P]$$ are dual.\end{theorem}

\begin{remark} As mentioned earlier, $M_1(X_{\bu})$ is not always a
1-motive: $N$ may contain torsion. But it is still possible to 
define the dual of $M_1(X_{\bu})$ (Deligne \cite{h} \S
10.2).\end{remark} 

\nd {\bf Proof.} It will be based on a number of lemmas. 
Let $B$ denote the maximal abelian quotient of $G$; we define the
torus $T$ by
means of the exact sequence 
$$ 0 \to T \to G \to B \to 0.$$
\begin{lemma} The abelian schemes $B$ and $Q$ are dual.
\end{lemma}

\nd {\bf Proof.}  The semiabelian scheme $G$ is defined by the top horizontal exact sequence
$$ 
\begin{CD}
A_{X_1}^0 @>d>> A_{X_0}^0 @>>> G @>>> 0\\ 
@V{j_1}VV @V{j_0}VV @VVV @. \\
A_{X'_1}^0 @>d>> A_{X'_0}^0 @>>> B @>>> 0;\\
\end{CD}
{} $$
the vertical maps are surjective, taking each semiabelian scheme to
its maximal abelian quotient.  
The abelian scheme  $Q$ is 
the neutral component of the (reduced) kernel of $Pic^{0, red}_{X'_0}
\xra{d_0-d_1}  Pic^{0, red}_{X'_1}$. The classical duality of the
Albanese and Picard varieties --- for any smooth
projective scheme $Z$, the abelian schemes $A^0_Z$ and $Pic^{0,red}_Z$
are naturally dual (compatible with morphisms). Applying this to  $X'_0$ and
$X'_1$, and using functoriality,  we see 
that $Q$ and $B$ are dual. \qed\ms

\begin{lemma}\label{tars}  The torus ${\mathcal Hom}(R, \gm)$ is
isomorphic to $T$.\end{lemma}

\nd {\bf Proof.} In the previous commutative diagram, the torus
corresponding to the kernel  of $j_1$ (resp. $j_0$) was identified to be
${\mathcal Hom}(R"_1,\gm)$ (resp. ${\mathcal Hom}(R"_0,\gm)$)  by
 Corollary \ref{242} (a) applied to $X_0$ and $X_1$. Consequently
 ${\mathcal Hom}(R,\gm)$ can  be identified with $T$.\qed\ms

\begin{lemma} There is a canonical isomorphism of tori:
$${\mathcal Hom} (N, \gm)  \cong C'. {} $$
\end{lemma}

\nd {\bf Proof.}  Fix a smooth proper scheme $W$ with projection $g:W \ra S$. 
The torus ${\mathcal Hom}(D_W,\gm)$ also 
represents the $fppf$ sheaf $R^0g_*\mathcal O^*_W$; a character of $D_W$
provides a nonzero function, constant (since $W$ is proper) on each
connected component of $W$, i.e., on  
each irreducible component since $W$ is smooth. The affine group
scheme $C'$ represents the homology sheaf of the 
complex :
\begin{equation}\label{4-7}
\mathcal R^0a_{0*}\mathcal O^*_{X'_0}  \xra{\delta^*_0 -
 \delta^*_1} \mathcal
R^0a_{1*}\mathcal O^*_{X'_1} \xra{\delta^*_0 -\delta^*_1 + \delta^*_2}
 \mathcal R^0a_{2*}\mathcal O^*_{X'_2}; \end{equation} this can be
 identified with the 
 following complex: 
\begin{equation}\label{4-8}
 {\mathcal Hom}(D_{X'_0},\gm) \xra{\delta^*_0 - \delta^*_1}
 {\mathcal Hom}(D_{X'_1},\gm)
 \xra{\delta^*_0 -\delta^*_1 + \delta^*_2}  {\mathcal
 Hom}(D_{X'_2},\gm). \end{equation} 
 Thus, $C'$ represents the homology sheaf of (\ref{4-8}).  We
 have defined $N$ to be the
 \'etale group  scheme representing the homology of the complex
 $$ D_{X_2} \xra{\delta_0 -\delta_1 + \delta_2} D_{X_1} \xra{\delta_0
 - \delta_1}  D_{X_0}. $$
 But note $D_{X_n} \cong D_{X'_n}$ since the complement $E_{\bu}$ of $X_{\bu}$
 in $X'_{\bu}$ is a  
 simplicial divisor. As the functor ${\mathcal Hom}(-,\gm)$
 establishes an anti-equivalence of 
the category of group 
schemes of
multiplicative type and \'etale group schemes, we get that $C' \cong
 {\mathcal Hom} (N,\gm)$.\qed\ms

A consequence is the isomorphism $Hom(C',\gm) \cong N$. 
In the next lemma, $\pi$ denotes the projection $P \ra Q$ (as above). 
 
\begin{lemma}\label{415} The 1-motive  $[R \xra{\pi~ h} Q] $ is dual to
 $[0 \ra G]$.\end{lemma}

\nd {\bf Proof.}  As observed earlier,  the 1-motive $[R \xra{\pi~ h} Q] $ 
 is the
$$ Ker([R''_0 \ra Pic^0_{X'_0}] \xra{\delta_0^* -\delta^*_1} [R''_1
\ra Pic^0_{X'_1}]),$$
rendering its dual to be the 1-motive $[0 \ra G]$, by
Theorem \ref{serth} applied to $X_0$ and $X_1$. \qed\ms

For any 1-motive $M$, Deligne (\cite{h} 10.2.11,
p. 67) defines the dual 1-motive $M^*$ by $M^*:= {\mathcal
RHom}(M,\gm[-1])$ where ${\mathcal RHom}$ is the 
derived functor of $Hom$ on the abelian category of complexes of group schemes
over $S$.
We follow his prescription in determining the dual of $M:=
M_1(X_{\bu})$:\ms

(a) The dual of $M/{W_{-2}(M)}$\ms

Recall that $P(V)= Pic^{0,red}_{X_{\bu}'}(V)$ classifies isomorphism classes of
pairs $(\mathcal L,
\alpha)$ with $\mathcal L$ an invertible sheaf on $X'_0\times_S V$ and
$\alpha$ an 
isomorphism $\delta^*_0(\mathcal L) \cong \delta^*_1(\mathcal L)$ over
$X'_1\times_S V$
satisfying a cocycle relation (\ref{rowss}) on $X'_2 \times_SV$. We may also interpret
$\alpha$ as a trivialization of the invertible sheaf $\delta^*_0(\mathcal
L) \otimes (\delta^*_1(\mathcal L))^{-1}$, i.e.,  an isomorphism $\delta^*_0(\mathcal
L) \otimes (\delta^*_1(\mathcal L))^{-1} \cong \mathcal O_{X_1 \times_S V}$ (also
required to satisfy a  cocycle condition).

There is a natural morphism of functors ${\mathcal Ext}^1([N
\xra{\partial} B], \gm) \xra{\psi} P$: elements of
the former are isomorphism classes of pairs $(\mathcal L, \alpha)$ where
$\mathcal L$ is an invertible sheaf on $X'_0$, defined by an element of $Q={\mathcal Ext}^1(B,\gm)$ (the dual
abelian scheme), and $\alpha$ is a trivialization of $\partial^*(\mathcal
L)$. The element $\mathcal L$ naturally extends to an invertible sheaf on
$X'_{\bu}$ because $\mathcal L \in Q$. Since $\partial$ is defined via the map 
$\delta_0-\delta_1: A_{X_1} \ra A_{X_0}$, we see that a
trivialization of $\partial^*(\mathcal L)$ amounts to a trivialization of
the invertible sheaf $\delta^*_0(\mathcal
L) \otimes (\delta^*_1(\mathcal L))^{-1}$ on $X'_1$. Comparing with
the definition of $P$ (recalled above), it is clear that the 
morphism $\psi$ is actually an isomorphism of functors. The 1-motives $[N
\xra{\partial} B]$ and $[0 \ra P]$ are dual.\ms

(b) For each $r\in R \cong Hom(T,\gm)$,
the extension $M$ of $M/{W_{-2}(M)}$ by $T$ defines an extension of
$M/{W_{-2}(M)}$ by $\gm$, and hence an element $s(r) \in P$. This
combines into a homomorphism $s: R \to P$. Therefore, the dual of
$M_1(X_{\bu})$ is the 1-motive $[R \xra{s} P]$.\ms

Thus, to finish the proof of Theorem \ref{duality}, it suffices to
show that $h=s$.\ms

Let $T_0$ be the torus
$Hom(R''_0,\gm)$. It is the kernel of the projection $j_*: A_{X_0} \ra
A_{X'_0}$ (see \ref{Serr}).   One may view
 $A^0_{X_0}$ as the total space of a $T_0$-bundle over $A^0_{X'_0}$. Similarly, $L$ (the universal torsor for
 $X'_0$) is the total space of a $T_0$-bundle over $L'$. Fix the universal morphisms $u: X_0 \ra L$,
 $u': X_0' \ra L'$  and
 the isomorphisms
\begin{equation}\label{haal}
Pic^{0,red}_{A^0_{X'_0}} \cong Pic^{0,red}_{L'} \cong
Pic^{0,red}_{X'_0}. {}
\end{equation}

 It is a complete tautology that both the $T_0$-bundles have the following property: For any $r \in R''_0
 \cong Hom (T_0,\gm)$, the line bundle obtained via $r$ over
 $A^0_{X'_0}$ and $L'$  
 corresponds to the element $[\mathcal O(r)] \in Pic^{0,red}_{X'_0}$
 under the isomorphisms in (\ref{haal}). 
Lemma \ref{415} and the fact that $R \hra R''_0$ together show 
that $\pi(s(r)) = \pi(h(r))$. 

Consider the restriction of the $T_0$-bundle $\tilde{L}:= (u')^*L$
over $X_0'$ to 
the open subscheme $X_0$. The universal
morphism $u: X_0 \to L$ is defined by lifting the map $u': X_0 \to
L'$ by means of a consistent choice 
of sections of the $T_0$-bundle $\tilde{L}$. For any $r \in R$, the map
$u$ determines a rational section $\tau_r$ of $\mathcal O(r)$ on
$X_0'$ whose divisor is $r$. On $X_0$, the section $\tau_r$ of
$\mathcal O(r)$ is nowhere
vanishing and regular. Similarly, we have a rational section $\tau_{-r}$ of the
invertible sheaf $\mathcal O(-r)$. One has a \emph{a priori} rational
section of 
the invertible sheaf $\mathcal L_1:= \delta^*_0(\mathcal 
O(r))\otimes \delta^*_1(\mathcal O(r))^{-1}$ on $X'_1$ given by
$\zeta_r:= \delta^*_0(\tau_r)\otimes \delta^*_1(\tau_{-r})$; this is because
$\mathcal O(r) ^{-1} = \mathcal O(-r)$. But $\mathcal L_1$ is the
trivial invertible sheaf $\mathcal O_{X_1}$ on $X_1$ and the section
$\zeta_r$ on $X_1$ is the identity section. 

The element $s(r)$  corresponds to a pair $(\pi(s(r)), \beta_r)$ where $\beta_r$ is a trivialization 
 of the invertible sheaf corresponding to $\pi(s(r))$ over $N$. As
 remarked earlier, the element $\pi(s(r))$ 
corresponds to the invertible sheaf $\mathcal O(r)$ on $X'_0$. Utilizing the definition of the map $\partial$, we see that the pull
back  $\partial^*(\mathcal O(r))$ is the invertible sheaf $\mathcal
 L_1$ on $X'_1$. Since $\beta_r$ is the pullback (via $\partial$) of
 the trivialization $\tau_r$ of $\mathcal O(r)$ on $X_0$, we obtain
 that  $\beta_r = \zeta_r$ on $X_1$. Therefore, $\beta_r$ also
 corresponds to the 
 identity trivialization of $\mathcal O_{X'_1}$. Comparing with the
 definition of $h$ (recalled below), we obtain that $h(r) = s(r)$. 

	One could alternatively proceed using the fact that $D_{X_i} =D_{X_i'}$
	for any $i$, and the sequences (\ref{4-7}) and (\ref{4-8}). 

 Let us recall the definition of the map $h$: Any
 element $r \in R$ satisfies $\delta_0^*(r) = \delta^*_1(r)= r_1$
 (equality of divisors on $X_1$). 
Therefore, there is a canonical isomorphism $\alpha_r: \delta_0^*(\mathcal O(r))\cong \delta^*_1(\mathcal
O(r))$,  corresponding to the identity in $Hom_{X'_1}(\mathcal O(r_1), \mathcal O(r_1))$. 
The invertible sheaf $\delta_0^*(\mathcal O(r))\otimes 
(\delta^*_1(\mathcal O(r)))^{-1}$ is canonically isomorphic to $\mathcal O_{X'_1}$; the element $\alpha_r$,
interpreted as a trivialization of the former, corresponds to the identity section of $\mathcal O_{X'_1}$. 
The map $h$ sends $r$ to the pair $(\mathcal O(r), \alpha_r)$.
 
    Consequently, the dual 1-motive of $M_1(X_{\bu})$ is $M^1(X_{\bu})$. \qed\ms 
{}

\begin{corollary} The Picard 1-motive $M^1(U_{\bu})$ is a
contravariant functor.\end{corollary}

\nd {\bf Proof.} Combine Theorem \ref{duality} with  the (evident)
 covariant functoriality of the 
 Albanese 1-motives. 

\begin{remark} Theorem \ref{duality} depends on the following two
 facts: 

(1) For any $r \in R$, the pullback $\delta^*_0(\mathcal O(r))\otimes
 \delta^*_1(\mathcal O(r)^{-1})$ is the trivial invertible sheaf
 $\mathcal O_{X'_1}$ on $X'_1$; this is by definition of the group
 $R$. 

(2) the pullback of any rational section $\tau_r$ of $\mathcal O(r)$
    (such that the divisor of $\tau_r$ is precisely $r$) under the map
    $\partial$ is the identity section of $\mathcal O_{X'_1}$; this is
    by definition of the map $\partial$.
\end{remark}
    
\section{Curves}\label{kurv}
A curve denotes a scheme of pure dimension one.
  For any curve $C$ over an algebraically closed field, there is the {\it
motivic $H^1$} of $C$ denoted $H^1_m(C)(1)(1)$ defined by Deligne
\cite{h} \S 10.3, and the 
homological 1-motive $h_1(C)$ defined by Lichtenbaum \cite{L1}; it is
  clear by inspection that  their definitions are
valid over a perfect field $k$. 

In this section, we provide a new proof (based on Theorem
\ref{duality}) of an unpublished theorem
of Lichtenbaum that 
$h_1(C)$ is dual to $H^1_m(C)(1)$.

We refer to (\cite{BLR} pp. 247-248; \cite{swan}; \cite{li3}) for the
details of seminormalization. We recall that the seminormalization
$\hat{C}$ of a curve $C$  is the
largest curve between the normalization $\tilde{C}$
and $C$ which is universally homeomorphic to $C$.  The association of
$\hat{C}$ and $\tilde{C}$  with $C$ is compatible with base change;
this requires $k$ to be perfect.

The  1-motives $H^1_m $ and $h_1$ of curves do not change upon
seminormalization and 
depend only on
the underlying reduced scheme structure. Therefore, we may assume that we deal
with reduced curves which are seminormal (i.e., that they do not have cuspidal
singularities).

\subsection{Jacobians of proper curves: a review}\label{prop}

Let $L$ be the (reduced) closed subscheme of a seminormal proper curve $C$ 
corresponding to the singular locus of ${C}$. Take $g:\tilde{C} \ra
 C$ to be the canonical projection. Put $K := L
\times_{C} \tilde{C}$; it is a reduced closed subscheme of
$\tilde{C}$. There is a natural
homomorphism  $g: D_K \ra D_L$ with kernel $D'$. The inclusion of $K$
in $\tilde{C}$ provides a homomorphism $D_K \ra D_{\tilde{C}}$. 
This yields a map $D' \xra{\mu}
D_{\tilde{C}}$. 

Since $\tilde{C}$ is normal, the
connected components of $\tilde{C}$ are the irreducible components of
$\tilde{C}$. So one has an isomorphism
$\pi_0({Pic}_{\tilde{C}}) \iso \pi_0(A_{\tilde{C}}) =
D_{\tilde{C}}$. Using this, one can state the Abel-Jacobi theorem as
follows: 

\begin{theorem} {\rm (Abel-Jacobi)} For any smooth projective curve
$X$, one has a natural isomorphism $$A_X \iso Pic_X;$$ this
isomorphism is induced by the map $X \hra Pic_X$ (each point of $X$ is
viewed as a divisor).\end{theorem} 

Consider the exact sequence of sheaves on ${C}$:
\begin{equation}\label{thak} 1 \ra \mathcal O^*_{{C}} \ra g_*\mathcal O^*_{\tilde{C}} \ra g_*\mathcal
O^*_{\tilde{C}}/{\mathcal O^*_{{C}}} \ra 1~.{} \end{equation}  
The quotient sheaf $\mathcal V:= g_*\mathcal
O^*_{\tilde{C}}/{\mathcal O^*_{{C}}}$ is supported on $L$.
Let $f: {C} \ra S$ be the structure morphism. The sequence
(\ref{thak}) 
 furnishes an exact sequence of sheaves on $S_{FL}$:
\begin{multline}\label{4-15} 
1 \ra R^0f_* \mathcal O^*_{{C}} \ra   R^0f_* (g_*\mathcal
O^*_{\tilde{C}})  \xra{\nu} R^0f_* \mathcal V \ra\\
 \ra R^1f_* \mathcal O^*_{\hat{C}} \ra   R^1f_* (g_*\mathcal
O^*_{\tilde{C}}) \ra 1
\end{multline}  
Since $g$ is finite, its higher derived functors are zero. Thus
$R^1f_* (g_*\mathcal O^*_{\tilde{C}}) = R^1(fg)_*\mathcal O^*_{\tilde{C}}$. In
other words, the last sheaf in (\ref{4-15}) corresponds to
the Picard scheme of $\tilde{C}$. 

The sheaf $R^0 f_* \mathcal V$ is the torus dual to the group  $D'$.
The sheaf $R^0f_* (g_*\mathcal O^*_{\tilde{C}})$ is the torus dual to
$D_{\tilde{C}}$. The map $\nu$ between these two sheaves in  corresponds
to the pullback of functions by the natural inclusion of $ K \hra
\tilde{C}$. Therefore $\nu$ is dual to the map $D' \xra{\mu}
D_{\tilde{C}}$ mentioned above. The cokernel of
$\nu$ is the torus $T$ dual to the group scheme  $D:= Ker(D' \xra{\mu}
D_{\tilde{C}})$. 

 An immediate corollary of (\ref{4-15})  is that
the Picard scheme of ${C}$ also
represents the functor which assigns, to a $S$-scheme $Z$, the group of
isomorphism classes of pairs ($\mathcal L, \alpha$) where $\mathcal L$ is an
invertible sheaf on $\tilde{C} \times_S Z$ and $\alpha$ is a
trivialization of the restriction of $\mathcal L$ to $K \times_S Z$. 
{}
\begin{remark} Note that the Picard scheme of a proper curve is always reduced
\cite{G1} Proposition 2.10. \end{remark} 

\subsection{1-motives of curves} 

Let $C$ be a curve (assumed to be seminormal).
Denote by $C'$ the unique
proper curve containing $C$ as a dense open subscheme such that the
singular locus of $C'$ is the same as that of $C$. As usual, 
$\tilde{C}$ and $\tilde{C'}$ denote their normalizations
\footnote{Note that there is no ambiguity here: $\tilde{C'} =
(\tilde{C})'$}. We may also  
describe the latter as the unique smooth proper curve containing
$\tilde{C}$ as an open dense subscheme. Denoting by $E$ (resp. $E_0$) the
closed complement of $C$ (resp $\tilde{C}$) in $C'$
(resp. $\tilde{C'}$), there is a natural projection $c: E_0 \ra E$, an
isomorphism of finite \'etale schemes. One has an isomorphism
\begin{equation}\label{cees}
c^*: {\ub{\Z}}^{w(E)} \iso
{\ub{\Z}}^{w(E_0)}.
\end{equation} 

Consider the natural morphism  $\alpha: \ub{\Z}^{w(E)} \ra
Pic{}_{C'} \ra \pi_0(Pic{}_{C'})$; let $v$ denote the induced map
$Ker(\alpha) \ra Pic^0_{C'}$.    
{}
\begin{definition} {\rm (Deligne)} {\it The motivic $H^1$ of $C$ is the 1-motive 
$$ H^1_m(C)(1):= [Ker(\alpha) \xra{v} Pic^0_{C'}]. {} $$
} \end{definition}

We recall that the special simplicial scheme
$ ((\tilde{C}/C)^{\Delta_n})_{n\geq 0}$ is smooth (the $(n+1)$-fold fibre
product $(\tilde{C}/C)^{\Delta_n}$ is the disjoint sum of
$\tilde{C}$, the diagonal, and a finite \'etale scheme) and that it
admits $((\tilde{C'}/{C'})^{\Delta_n})_{n \geq 0}$ as a smooth
compactification with a (smooth!) simplicial
divisor $E_{\bu}$ as complement. Let us put $$ C_n:= (\tilde{C}/C)^{\Delta_n}, \qquad
C'_n:= (\tilde{C'}/{C'})^{\Delta_n} \qquad (n \geq 0). {} $$
Both of the schemes $C_{\bu}$ and $C'_{\bu}$ are smooth simplicial schemes. We have $C_0= \tilde{C}$
and $C'_0 = \tilde{C'}$.  

The natural augmentations $a:C_{\bu} \ra C$ and $b:C'_{\bu}
\ra C'$ are proper hypercoverings of $C$ and $C'$ respectively.

{}
\begin{lemma}\label{rosen} The map $b^*:Pic^0_{C'} \ra
Pic^0_{C'_{\bu}}$ is an isomorphism.
\end{lemma}

\nd {\bf Proof.} The map from $b_0: \tilde{C'}= C'_0 \ra C'$ is
finite and surjective. Interpreting the conditions on the elements of
$Pic^0_{C'_{\bu}}$ as descent data, we see that the map $b^*$
is injective (\cite{BLR}, \S6) (see also \ref{prop}).
 
 For the surjectivity, one uses the well known fact that descent
data for invertible sheaves is effective for finite surjective
morphisms (\cite{BLR}, \S 6.5, Theorem 1, p. 157).\qed\ms

Let $F$ denote the singular locus of $C$; we take $F_0:= F \times_C
\tilde{C}$.  
Denote by $V$ the proper seminormal curve corresponding to the
modulus $E_0$ on $\tilde{C}'$ (\cite{Se1} IV, no. 4, p. 61 and
76). There is a distinguished $S$-point $U$
 which is the unique singular point of $V$. The pullback of $U$ under the map
$\tilde{C}' \ra V$ is $E_0$. The curve $\tilde{C}'$ is the
normalization of $V$. We set $H:= Pic^0_V$.

We recall the following
{}
\begin{theorem}\label{rose} {\rm (Rosenlicht-Serre)} There is a natural
isomorphism $$H \iso A^0_{\tilde{C}} {} $$ of semiabelian schemes
over $S$.\end{theorem}

\nd {\bf Proof.} This follows by combining\ms 

 $\bu$ Theorem 1 in (\cite{Se1}
Chapter V, p.88):  $H$ is isomorphic to the
generalized Jacobian of $\tilde{C}'$ corresponding to the modulus $E_0$
and\ms

 $\bu$ Example 1 in \cite{Se3}:  the generalized
Jacobian of $\tilde{C}'$ corresponding to the modulus $E_0$ is none
other than $A^0_{\tilde{C}}$. \qed\ms 

Notice that  $F$ and $E$ (similarly $E_0$ and $F_0$) are
disjoint. 

Put $D':= Ker(D_{F_0} \ra D_F)$; there is a natural map $D'
\xra{\mu} A_{\tilde{C}}$
induced by the inclusion of $F_0$ in $\tilde{C}$. We set $D:= Ker (D'
\xra{\mu} A_{\tilde{C}} \ra D_{\tilde{C}})$, the last being the group of
connected (= irreducible) components of $\tilde{C}$. The group $D(\bar{S})$ is
 generated by differences $(x-y)\in D_{F_0}(\bar{S})$ such that
 $b_0(x) = b_0(y) \in F(\bar{S})$ and $x,y$ are in the same
 geometric component of $\tilde{C}$.

Using Theorem \ref{rose}, we may define $h_1(C)$ \cite{L1}.  

\begin{definition} {\rm (Lichtenbaum)} The homological 1-motive
$h_1(C)$ is the 1-motive
$$ [D \xra{\mu} A^0_{\tilde{C}}]. {} $$ \end{definition}

The key identification is in the
\begin{theorem}\label{curv}  One has natural isomorphisms 

{\rm (i)} $ H^1_m(C)(1) \iso  M^1(C_{\bu}).$ 
 
{\rm  (ii)} $ h_1(C) \iso M_1(C_{\bu}).$ 
\end{theorem} 

\nd {\bf Proof.} We begin with (i). 
We have defined $R$ to be the subgroup scheme  of $\ub{\Z}^{w(E_0)}$
corresponding to the intersection of the kernels of 
$\ub{\Z}^{w(E_0)} \ra \pi_0(Pic^{0}_{C'_0})$ and
$ \ub{\Z}^{w(E_0)} \xra{\delta^*_0 -\delta^*_1} \ub{\Z}^{w(E_1)}$.

It is a simple consequence of the definition of $C'_{\bu}$ that we
have a complex:
$$ {\ub{\Z}}^{w(E)} \xra{c^*} {\ub{\Z}}^{w(E_0)} \xra{\delta^*_0 -\delta^*_1}
\ub{\Z}^{w(E_1)}. {} $$ But $c^*$ is an isomorphism (\ref{cees}). Therefore, $R= Ker
(\ub{\Z}^{w(E_0)} \ra \pi_0(Pic^{0}_{C'_0}))$.  Since the group
of connected components of the Picard schemes of $C'$ and $C'_0
=\tilde{C'}$ are isomorphic (\cite{BLR}, Proposition 10, \S9.2,
p.248-9), we find  an isomorphism $c^*: Ker(\alpha) \iso R$.
 
The diagram below is commutative:
$$
\begin{CD}
Ker(\alpha) @>{v}>> Pic^0_{C'}\\
@V{c^*}VV @V{b^*}VV\\
R @>{h}>> Pic^0_{C'_{\bu}}\\
\end{CD}
{} 
$$
with isomorphisms (see Lemma \ref{rosen}) for vertical arrows. This
proves (i).  

We remark that  $C_0:= \tilde{C}$
and $C_1:= \tilde{C}\times_C \tilde{C}$. As usual, we write the 
Albanese 1-motive of $C_{\bu}$ as $[N \xra{\partial} G]$ (see \S
\ref{ilbee}). 

The Albanese 1-motive of a special scheme
$W$ is the neutral 
component of the Albanese scheme $A_W$ of $W$ (\ref{Serr}). Any normal
curve is a 
special scheme. 
It is clear from the definitions of the Albanese 1-motive and
the homological 1-motive that, for a normal curve $Y$, one has $h_1(Y) =
M_1(Y) = A^0_{Y}$.  

Returning back to our curve $C$ and $C_{\bu}$, let us 
set $A_i:= A_{C_i}$. Since each $C_i$ is the disjoint union of the
diagonal (isomorphic to $\tilde{C}$) and a finite \'etale scheme $X_i$
($i \geq 1$)\footnote{This is because $C_0 = \tilde{C}$.}, we
have that $A_i= A_{\tilde{C}} \times_S D_{X_i}$. We have $A_0= A_{\tilde{C}}$.
Therefore, the
simplicial large group $A_{\bu}$ contains a constant simplicial large
group $B_i:= A_{\tilde{C}}$ (for all non-negative integers
$i$). Viewing these simplicial groups (and $D_{X_i}$) as complexes (by normalization),
we see 

\nd (1) that  $A_{\tilde{C}}= H_0(B_{\bu}) = H_0(A_{\bu})$. 

\nd (since $C_0 = \tilde{C}$) that  $\pi_0(A_0) =
D_{\tilde{C}}$. This in turn implies that 
$$ N = \frac{Ker~(D_{X_1} \xra{\delta_0 -\delta_1}
D_{\tilde{C}})}{Im~(D_{X_2} \xra{\delta_0 -\delta_1 + \delta_2}
D_{X_1})}.$$ 

In particular, we deduce $G = A^0_{\tilde{C}}$. 
So we have already proved that $W_{-1}h^1(C) = W_{-1}M^1(C_{\bu})$. We will use the
previous description of $N$ to show that 
$D$ is actually $N$ (compatible with the map $\mu$) and this will
finish the proof of (ii). 

Let us put $F_n$ to be the $(n+1)$-fold fibre product of $F_0$ over
$F$. We have $F_1:=F_0 \times_F F_0$ and $F_2:= F_0 \times_F F_0
\times_F F_0$. 
 
The group $D_{X_1}(\bar{S})$ is generated by points $(x,y) \in
 F_1(\bar{S})$ with $x \neq y$ (i.e. $x,y$ are distinct points of $F_0$). By definition, the points $x, y \in
 \tilde{C}(\bar{S})$ satisfy $b_0(x)=b_0(y) \in C(\bar{S})$. Set
 $K:= Ker(D_{X_1} \xra{\delta_0 -\delta_1} D_{\tilde{C}})$. The group
 $K(\bar{S})$ is generated by points $(x,y) \in F_1(\bar{S})$ which
 map to the same (geometrically) connected component of $\tilde{C}$.
There is a natural morphism $\gamma: K \to D$, which we describe on
 the level of the geometric points: we send $(x,y)$ to the element
 $x-y\in D(\bar{S})$. The morphism $\gamma$ is surjective (by the
 seminormality of $C$). 

Let us turn
 to the group $D_{X_2}(\bar{S})$. It is generated by elements $(x,y,z)
 \in F_2(\bar{S})$ with $x\neq y$, $y\neq z$, $x \neq z$. By
 definition, we have $b_0(x)=b_0(y)=b_0(z) \in C(\bar{S})$. One has a
 map
$$\delta: D_{X_2} \xra{\delta_0 -\delta_1 + \delta_2} K \qquad (x,y,z)
 \mapsto (y,z) - (x,z) + (x,y). {} $$ 
The composition $\gamma . \delta$ can be checked to be zero; in fact,
 $Ker(\gamma)= Im(\delta)$. Since the group $N$ is $Coker(\delta)$,
 we see that $N = D$.
It is clear that the maps from $N$ and $D$ to $ A^0_{\tilde{C}}$ are
 compatible with this isomorphism. This proves (ii). \qed\ms {}

\begin{theorem}\label{slick} {\rm (Lichtenbaum, unpublished)} The 1-motives
  $H^1_m(C)(1)$ and $h_1(C)$ are dual. \end{theorem}

\nd {\bf Proof.} This follows from Theorems \ref{curv} and
\ref{duality}. \qed\ms

\begin{corollary} {\rm (Hodge Realization)} For a curve $C$
over $\C$, one has an isomorphism
$$H_1(C) \cong \mathfrak T(h_1(C)).$$ \end{corollary}
\nd {\bf Proof.} Combine Theorems 10.3.8 of
\cite{h} and \ref{duality}.   \qed\ms {}

The motivic $H^1$ of $C$ can also be interpreted as the motivic
relative cohomology group  
$H^1_m(\tilde{C}' - E_0$ rel $F_0)(1)$. Likewise, one can interpret
$h_1(C)$ as the motivic relative homology group  
$H_1^m(\tilde{C}' - F_0$ rel $E_0)$. With these identifications, the
duality theorem translates to a  Lefschetz 
duality theorem $H^1_m(\tilde{C}' - E_0$ rel $F_0)(1) \iso
(H_1^m(\tilde{C}' - F_0$ rel $E_0))^*$ (see 1.5 of \cite{dmot}).   
 The resulting duality of the Hodge realizations is a classical result
\cite{Ch1}, \cite{weilc} (especially p. 390-92).

\section{1-motives associated with varieties}\label{positive}
Let $Y$ be a  scheme over $S$. By de Jong \cite{dj}, there exists a diagram 
   
 $$
\begin{CD}
X_{\bu}      @>j>> X'_{\bu}  @<<< E_{\bu}\\
@V{\alpha}VV       @.            @. \\
Y            @.     {}      @.   {}\\
\end{CD}
{} $$
where the smooth simplicial scheme $X_{\bu}$ is  the open complement
of a reduced simplicial divisor with normal crossings $E_{\bu}$ in a
smooth projective simplicial scheme $X'_{\bu}$ such that, via $\alpha$, 
$X_{\bu}$ is a proper hypercovering of $Y$. Using this, it is possible
(as in \cite{h} \S 6.3)
to translate the theory developed for simplicial schemes such as
$X_{\bu}$ to obtain invariants for the scheme $V$. This has been
treated  in \cite{me2}.

\end{document}